\documentclass[10pt]{amsart}

\usepackage{a4,amssymb,mathrsfs,amscd}
\usepackage[arrow,matrix,curve]{xy}\SilentMatrices
\newdir{ >}{{}*!/-5pt/\dir{>}}

\DeclareMathAlphabet{\bsf}{OT1}{cmss}{bx}{n}

\def\Nil{{\bsf{Nil}}}
\def\nil{{\bsf{nil}}}

\newcommand{\qu}{{\sf Quad}}

\def\Ga{\Gamma}

\newcommand{\Ker}{{\sf Ker}}
\DeclareMathOperator{\Id}{{\sf Id}}
\DeclareMathOperator{\cok}{{\sf Coker}}
\DeclareMathOperator{\im}{{\sf Im}}

\def\Z{{\mathbb Z}}

\newcommand{\tp}{\otimes}

\def\ZZ{\mathsf{Z}}
\newcommand{\modr}{\mathsf{Mod_{\RR}}}
\newcommand{\cp}{\mathsf{CP}}
\newcommand{\cpr}{\mathsf{CP_{\RR}}}

\def\Ga{\Gamma}

\def\La{\Lambda}

\def\ga{\gamma}

\def\theodef#1#2{\newtheorem{#1}[subsubsection]{#2}}

\theodef{Lemma}{Lemme} \theodef{Corollary}{Corollaire}
\theodef{Proposition}{Proposition} \theodef{Theorem}{Th\'eor\`eme}
\theoremstyle{definition} 
\theodef{Definition}{D\'efinition}
\theodef{Remark}{Remarque} \theodef{Example}{Exemple}
\theodef{Examples}{Exemples}
\theodef{Exercise}{Exercice} \theodef{Exercises}{Exercices}

\newenvironment{MYitemize}{\begin{itemize}}{
\end{itemize}}

\numberwithin{equation}{subsection}

\newif{\ifmargnote}\margnotetrue
\newcommand{\marg}[1]%
{\ifmargnote\ifodd\value{page}\normalmarginpar\marginpar{\vtop{
\footnotesize #1}}
\else\reversemarginpar\marginpar{\vtop{\footnotesize #1}}\fi\else\relax\fi}

\newcommand{\MA}{\underline{M}}
\newcommand{\NB}{\underline{N}}

\newcommand{\RR}{\underline{R}}
\newcommand{\ud}[1]{\underline{#1}}
\newcommand{\annc}[1]{(\xymatrix{#1_e\ar[r]^H&#1_{ee}\ar[r]^P&#1_e})}
\newcommand{\Znil}{\Z_{\nil}}
\newcommand{\ovl}[1]{\overline{#1}}
\newcommand{\ti}{\text{-}}
\newenvironment{preuve}{\begin{proof}[Preuve]}{\end{proof}}

\begin{document}
\title{Morphismes quadratiques entre modules sur un anneau carr\'e}
\author[H. GAUDIER]{Henri GAUDIER}
\author[M. HARTL]{Manfred HARTL}
\address{Univ. LILLE Nord de France, F-59000 LILLE, France
\hfil\break\strut\hskip4.3mm
UVHC,LAMAV et FR CNRS 2956, F-59313 VALENCIENNES, France}
 \email{manfred.hartl@univ-valenciennes.fr}
\email{henri.gaudier@univ-valenciennes.fr}
\subjclass[2000]{18D50, 13C99, 20N99, 20F18}
\keywords{quadratic map, square ring, modules over square rings, commutative square ring, operad}
\date{\today}
\maketitle
\begin{abstract}
We introduce the notions of a commutative square ring $ \RR$ and of a quadratic map  between modules over $\RR$, called $\RR$-quadratic map. This notion generalizes various notions of quadratic maps between algebraic objects in the literature. We construct a category of quadratic maps between $\RR$-modules and show that it is a right-quadratic category and has an internal Hom-functor.

Along our way, we recall the notions of a general square ring $\RR$ and of a module over $\RR$, and discuss their elementary properties in some detail, adopting an operadic point of view. In particular, it turns out that the associated graded object of a square ring $\RR$ is a nilpotent operad of class $2$, and the associated graded object of an $\RR$-module is an algebra over this operad, in a functorial way. This generalizes the well-known relation between groups and graded Lie algebras (in the case of nilpotency class $2$). We also generalize some elementary notions from group theory to modules  over square rings.

\end{abstract}

\section{Introduction}
In this paper we generalize various notions of quadratic map between algebraic objects studied in the literature: quadratic forms certainly are the most classical example, followed by the notion of quadratic map between modules defined to be the sum of  a linear and a homogenous  quadratic map. The latter notion does not cover the map $\mathbb Z \to \mathbb Z$, $n\mapsto {n\choose 2}$, which, however, is a quadratic map in Passi's sense who first introduced polynomial maps from  groups to abelian groups  \cite{Pa68}. In \cite{GH} we introduce a more general notion of quadratic map between modules which in the case of $\mathbb Z$-modules is equivalent with Passi's notion specialized to abelian domains. On the other hand, quadratic maps between arbitrary groups were studied in \cite{Niq}, \cite{QmG} after they had arosen in the new field of ``quadratic algebra'', namely the theory of square rings and their modules, see \cite{BHP} and numerous papers by Baues, Jibladze, Muro and Pirashvili giving applications of these structures in homotopy theory, see for example \cite{B}, \cite{BHP}, \cite{BJP1}, \cite{BJP2}, \cite{BM}.  

Square rings and their modules also constitute the framework of this paper; this allows to define a notion of quadratic map which covers the one for modules over an ordinary ring in \cite{GH}, as well as quadratic maps between nilpotent groups of class $2$ (which is the generic case of quadratic maps between groups, see proposition \ref{quadNil2gen} below); but it also provides a notion of quadratic map between nilpotent algebras of class $2$ over any operad. We start by recalling the definition and elementary properties of square rings and their modules, noting in addition that the associated graded group of a square ring $\RR$ has the structure of a nilpotent  {operad} of class $2$, and that the 
associated graded group of an $\RR$-module is an algebra over this operad, in a functorial way. We need, however, to distinguish two notions of module over $\RR$, the original one from \cite{BHP}, and a slightly more general one coming from the theory of modules over a square ringoid, by specializing to the one-object case \cite{chap2}. The latter notion is needed in the last section section where we construct a category of 
$\RR$-quadratic maps between $\RR$-modules, equipped with an internal Hom-functor, thus generalizing the category $\mathsf{CP}$ of quadratic maps between pairs of groups in \cite{QmG}. The construction of this category indeed is the main objective of this paper: it allows us in forthcoming work to introduce square algebras and their modules, as a new approach to the study of quadratic maps, even between vector spaces.

To define  $\RR$-quadratic maps, $\RR$ has to be commutative, a property we introduce in section \ref{sec:AnnCarComm}, and which can be conveniently expressed in terms of the operad associated with $\RR$. Several equivalent characterizations of $\RR$-quadratic maps are given, as well as some canonical examples; as a striking fact, the endomorphism of $\RR$ given by  $r\mapsto r^2$ is no longer $\RR$-quadratic in general, nor is the map $r\mapsto r^2-r$, which is the universal quadratic map annihilating $1$ on an ordinary commutative ring. In fact, lack of a canonical $\RR$-quadratic endomorphism of the free $\RR$-module of rank $1$, we were not able so far to determine the structure of the target of the universal $\RR$-quadratic map on an $\RR$-module, as we did in \cite{GH} for quadratic maps on a module over an ordinary ring.

\section{Anneaux carr\'es et leurs modules}
\label{sec:AnnCarMod}
\bigskip
\subsection{Anneaux carr\'es}
\label{sec:AnnCar}{\strut}

\begin{Definition}\cite{BHP} \label{df:acarr}
Un \emph{anneau carr\'e} $\RR=\annc R$ est donn\'e par
\begin{MYitemize}
 \item un pr\'e-anneau (near-ring) \`a droite $R_e$ (i.e. un objet qui
satisfait tous les axiomes d'un anneau sauf la commutativit\'e de
l'addition et la distributivit\'e \`a gauche),
 \item un groupe ab\'elien $R_{ee}$,
\item une application quadriadditive(\footnote{\ Pour \'eviter toute
ambiguit\'e, une application qui respecte l'addition, sera dite additive
(et non pas lin\'eaire).}) $R_e\times R_e\times R_{ee}\times R_e \to
R_{ee}$, $(r,s,x,t)\mapsto (r,s)\cdot x\cdot t$ munissant $R_{ee}$ d'une
action \`a gauche du mono\"ide multiplicatif $R_e\times R_e$ et d'une
action \`a droite du mono\"ide multiplicatif $R_e$, compatibles entre
elles,(\footnote{\ ou si l'on pr\'ef\`ere d'une action \`a gauche du
mono\"ide $R_e\times R_e \times R_e^{op}$.})
\item un morphime de groupes $T:R_{ee}\to R_{ee}$ tel que $T^2=\Id$
\item une application $H:R_e\to R_{ee}$ et un morphisme de groupe
$P:R_{ee}\to R_e$ v\'erifiant les propri\'et\'es suivantes, o\`u $r,s,t\in
R_e$ et $x,y\in R_{ee}$:
\end{MYitemize}
\begin{align}
 \tag{AC0}\label{ax:AC0} &PHP=P+P,\\
\tag{AC1}\label{ax:AC1} &P\text{ est binaturel, i.e. }
P((r,r)\cdot x\cdot s)=rP(x)s,\\
\tag{AC2}\label{ax:AC2} &T=HP-\Id,\\
\tag{AC3}\label{ax:AC3} &PT=P,\\
\tag{AC4}\label{ax:AC4} &(P(x),r)\cdot y=(r,P(x))\cdot y=0=y\cdot
P(x),(\footnotemark)\\
\tag{AC5}\label{ax:AC5} &H(r+s)=H(r)+H(s)+(s,r)\cdot H(2),\text{ o\`u }
2=1_R+1_R,\kern55mm\\
\tag{AC6}\label{ax:AC6} &H(rs)=(r,r)\cdot H(s)+H(r)\cdot s,\\
\tag{AC7}\label{ax:AC7} &(r+s)t=rt+st+P((r,s)\cdot H(t)),\\
\tag{AC8}\label{ax:AC8} &T((r,s)\cdot x\cdot t)=(s,r)\cdot T(x)\cdot t.
\end{align} 
\footnotetext{\ la derni\`ere relation ne figure pas dans \cite{BHP}.}
\end{Definition}
 

 \begin{Remark}\label{Op(R)}
 Rappelons \cite{BHP} qu'alors $\ovl{R}=\cok P$ est un anneau et que $R_{ee}$ est un $\ovl{R}\otimes \ovl{R} \otimes \ovl{R}^{op}$-module. Gr\^ace \`a la  relation (AC8), on obtient ainsi une \textit{op\'erade nilpotente de classe $\le 2$}, not\'ee  OP$(\RR)$, o\`u  OP$(\RR)_1=\ovl{R}$ et   OP$(\RR)_2 = R_{ee}$ (rappelons qu'une op\'erade $\mathbb P$ est dite nilpotente de classe $\le n$ si $\mathbb P_k=0$ pour $k>n$).

\end{Remark}

Remarquons \'egalement que $H(1)=0$ (prendre $r=s=1$ dans (\ref{ax:AC6})), et que (\ref{ax:AC2}) $\Rightarrow$ ( (\ref{ax:AC0})
$\Leftrightarrow$ (\ref{ax:AC3}) ).

Des exemples d'anneaux carr\'es et de leur modules (cf.\ le paragraphe suivant) seront donn\'es dans le paragraphe \ref{sec:exBHP} plus loin.


\subsection{Modules sur un anneau carr\'e}
\label{sec:ModAnnCar}
Nous donnons, en parall\`ele, deux d\'efinitions pour la notion
de module sur un anneau carr\'e. La premi\`ere se trouve dans la
litt\'erature, par exemple dans \cite{BHP}, nous la notons ici BHP-module.
La seconde, inspir\'ee de la notion analogue pour les groupes (\cite{QmG}),
s'impose lorsque l'on veut que la somme et la compos\'ee d'applications
quadratiques restent quadratiques. C'est aussi celle qu'on obtient en
sp\'ecialisant la notion de module sur un ann\'elide carr\'e (\cite{B2}).
Nous la notons ici CP-module.

\smallbreak
\noindent\begin{minipage}[t]{0.48\textwidth}
\begin{Definition} (\cite{BHP} def 7.9)\label{df:BHPmodc}
 Soit $\RR=\annc{R}$ un anneau carr\'e. Un \emph{$\RR$-BHP-module \`a
droite} est un groupe $M$ (additif mais non n\'ecessaire\-ment commutatif)
muni d'op\'erations

\begin{MYitemize}
 \item[\llap{--\qquad}] $\kern-7mm M\times R_e\to M$,\qquad $(m,r)\mapsto
m\cdot r$,
\item[\llap{--\qquad}] $\kern-7mm M\times M\times R_{ee}\to M$,\qquad
$(m,n,x)\mapsto [m,n]\cdot x$ 
\end{MYitemize}
satisfaisant les relations suivantes, o\`u $r,s\in R_e$, $x,y\in R_{ee}$,
$m,m',n\in M$ :

\kern1.0mm
\begin{align}
\tag{MC1}\label{ax:MC1} &\kern-1mm m\cdot 1=m,\qquad (m\cdot r)\cdot
s=m\cdot (rs),\\
\notag &\kern-1mm  m\cdot(r+s)=m\cdot r+m\cdot s,\\
\tag{MC2}\label{ax:MC2} &\kern-1mm (m+n)\cdot r=m\cdot r+n\cdot
r+[m,n]\cdot H(r),\\
\tag{MC3}\label{ax:MC3} &\kern-1mm m\cdot P(x)=[m,m]\cdot x,\\
\tag{MC4}\label{ax:MC4} &\kern-1mm [m,n]\cdot T(x)=[n,m]\cdot x,\\
\tag{MC5}\label{ax:MC5} &\kern-1mm [m,n]\cdot x\text{ est additif en }m,\
n\text{ et }x,\\
\tag{MC6}\label{ax:MC6} &\kern-1mm([m\cdot r,n\cdot s]\cdot x)\cdot
t=[m,n]\cdot((r,s)\cdot x\cdot t),\\ 
\tag{MC7}\label{ax:MC7} &\kern-1mm [[m,m']\cdot x,n]\cdot y=0.
\end{align}
\end{Definition}
\end{minipage}\qquad 
\begin{minipage}[t]{0.48\textwidth}
\begin{Definition}(\cite{chap2})\label{df:CPmodc}
 Soit $\RR=\annc{R}$ un anneau carr\'e. Un \emph{$\RR$-CP-module \`a
droite} est un couple $(M,A)$ form\'e d'un groupe $M$ et un sous
groupe $A$ de $M$, muni d'op\'erations
\begin{MYitemize}
 \item[\llap{--\qquad}] $\kern-7mm M\times R_e\to M$,\qquad $(m,r)\mapsto
m\cdot r$,
\item[\llap{--\qquad}] $\kern-7mm M\times M\times R_{ee}\to M$,\qquad
$(m,n,x)\mapsto [m,n]\cdot x$ 
\end{MYitemize}
satisfaisant les relations suivantes, o\`u $r,s\in R_e$, $x,y\in R_{ee}$,
$m,m',n\in M$ :
\begin{align}
\tag{MC0}\label{ax:MC0} &\kern-1mm A\text{ est stable sous l'action de
}R_e\\
\tag{MC1} &\kern-1mm m\cdot 1=m,\qquad (m\cdot r)\cdot s=m\cdot (rs),\\
\notag & \kern-1mm m\cdot(r+s)=m\cdot r+m\cdot s,\\
\tag{MC2} &\kern-1mm (m+n)\cdot r=m\cdot r+n\cdot r+[m,n]\cdot H(r),\\
\tag{MC3} &\kern-1mm m\cdot P(x)=[m,m]\cdot x,\\
\tag{MC4} &\kern-1mm [m,n]\cdot T(x)=[n,m]\cdot x,\\
\tag{MC5} &\kern-1mm [m,n]\cdot x\text{ est additif en }m, n\text{ et }x,\\
\tag{MC6} &\kern-1mm ([m\cdot r,n\cdot s]\cdot x)\cdot
t=[m,n]\cdot((r,s)\cdot x\cdot t), 
\end{align}
\vglue-6.1mm
\begin{align}
\tag{MC7a}\label{ax:MC7a} &\kern-12mm 
[m,n]\cdot x=0\text{ si }m\in A,\\
\tag{MC7b}\label{ax:MC7b} &\kern-12mm  [m,n]\cdot x\in A.
\end{align}
\end{Definition}
\end{minipage}
\smallbreak

\begin{Remark}
 L'axiome \ref{ax:MC7} d\'ecoule clairement des axiomes [\ref{ax:MC7a}] et
[\ref{ax:MC7b}], par cons\'equent tout $\RR$-CP-module est un
$\RR$-BHP-module.
\end{Remark}

\begin{Remark}
 La d\'efinition peut s'interpr\'eter de la fa\c con suivante. On se donne
un groupe $M$ (dont on verra plus tard (\ref{pr:modelem}) par l'axiome
\ref{ax:MC7} qu'il est nilpotent de classe $\leq 2$) muni d'une famille
d'op\'erations unaires param\'etr\'ees par $R_e$ et d'une famille
d'op\'erations binaires param\'etr\'ees par $R_{ee}$ (on verra plus loin 
que quand $\RR$ est commutatif les op\'erations unaires sont quadratiques, 
et les op\'erations binaires sont bilin\'eaires). Les applications
$H$ et $P$ relient ces deux familles d'op\'erations: l'axiome \ref{ax:MC2}
exprime que l'effet crois\'e de l'op\'eration unaire associ\'ee \`a 
$r\in R_e$ est l'op\'eration binaire associ\'ee \`a $H(r)$, et l'axiome
\ref{ax:MC3} exprime que la compos\'ee de l'application diagonale de $M$ et
de l'op\'eration binaire associ\'ee \`a $x\in R_{ee}$ est l'op\'eration
unaire associ\'ee \`a $P(x)$.
\end{Remark}

\begin{Examples}\label{ex:RmodRe}\ 

\noindent\begin{minipage}[t]{0.48\textwidth}
 1)\ On v\'erifie ais\'ement que le groupe $R_e$ muni des actions
\begin{align*}
 R_e\times R_e&\to R_e& r\cdot s&=rs\\
R_e\times R_e\times R_{ee}&\to R_e& [r,s]\cdot x&=P((r,s)\cdot x),
\end{align*}
est un $\RR$-BHP-module.
\end{minipage}
\quad 
\begin{minipage}[t]{0.48\textwidth}
 2)\ Soit $A$ un sous groupe de $R_e$ stable par la multiplication \`a
droite par $R_e$ et tel que $P(R_{ee})\subseteq A\subseteq \ZZ_{\RR}(R_e)$
(\footnote{\ cette notation est d\'efinie en \ref{df:Rcentre}.}). On
v\'erifie ais\'ement que $(R_e,A)$ muni des m\^emes actions que ci-contre
est un $\RR$-CP-module.
\end{minipage}
\smallbreak

3)\ Soit $M$ un $\ovl R$-module \`a droite, par restriction des scalaires
on a une action \`a droite de $R_e$ sur $M$. En prenant l'application nulle
$M\times M \times R_{ee}\to M$, on obtient une structure de
$\RR$-BHP-module sur $M$. On v\'erifie ais\'ement que pour tout sous-$\ovl
R$-module \`a droite $M'$, le couple $(M,M')$ est un $\RR$-CP-module. En
particulier, si $\RR$ est un anneau carr\'e, $R_{ee}$ est naturellement un
$\RR$-BHP-module. Et pour tout $R_e$-sous-module $B$ de $R_{ee}$,
$(R_{ee},B)$ est un $\RR$-CP-module.
\end{Examples}
D'autres exemples seront donn\'es plus loin en \ref{sec:exBHP}.

\begin{Proposition}\label{pr:modelem}{\strut }
\emph{\cite{chap2} } Soit $M$ un $\RR$-BHP-module, $m,n\in M$,
$r\in R_e$ et $x,y\in R_{ee}$ :
 \begin{enumerate}
   \item les \'el\'ements $[m,n]\cdot x$ sont dans le centre de $M$,
   \item on a $[m,n]\cdot H(2)=[n,m]$,
   \item le groupe $M$ est nilpotent de classe 2,
   \item $(-m)\cdot r=-(m\cdot r)+[m,m]\cdot H(r)=-(m\cdot r)+m\cdot
PH(r)$,
   \item $([m,n]\cdot x)\cdot P(y)=0$, en particulier $P(x)P(y)=0$.
 \end{enumerate}
\end{Proposition}
On peut interpr\'eter la propri\'et\'e 2), en disant que les op\'erations
binaires associ\'ees aux \'el\'ements $x$ de $R_{ee}$ g\'en\'eralisent les
commutateurs de $M$ (pour l'addition) d\'efinis par $[m,n]:=m+n-m-n$.
\begin{preuve} Remarquons d'abord que l'axiome \ref{ax:MC2} donne
\begin{equation*}
[m,n]\cdot H(2)=-n-n-m-m+m+n+m+n=-n+[-n,-m]+n
\end{equation*}
Rempla\c cant $m$ et $n$ par leurs oppos\'es on a gr\^ace \`a
(\ref{ax:MC5})
\begin{equation}
 [m,n]\cdot H(2)=[-m,-n]\cdot H(2)=n+[n,m]-n\tag{*}
\end{equation} 
L'axiome \ref{ax:MC7} avec $y=H(2)$ et la relation (*) donnent alors
\begin{align*}
 0&=[[m,m']\cdot x,n]\cdot H(2)\\
&=n+[n,[m,m']\cdot x]-n
=[n,[m,m']\cdot x],
\end{align*}
ce qui prouve (1).

Puisque $[m,n]\cdot H(2)$ est dans le centre, la propri\'et\'e (2)
d\'ecoule imm\'ediatement de (*).
Puisque les commutateurs sont centraux, la propri\'et\'e (3) est vraie.

Pour la quatri\`ere relation, on a :
\begin{equation*}
 0=(m-m)\cdot r=m\cdot r+(-m)\cdot r +[m,-m]\cdot H(r) = m\cdot r+(-m)\cdot
r -[m,m]\cdot H(r),
\end{equation*}
ce qui montre la formule, en terminant avec (\ref{ax:MC3}).
Pour la cinqui\`ere relation, on a d'apr\`es la derni\`ere relation de
\ref{ax:AC4} 
\begin{equation*}
 ([m,n]\cdot x)\cdot P(y)=[m,n]\cdot (x\cdot P(y))=[m,n]\cdot 0=0.
\end{equation*}
En prenant alors $M=R_e$, et $m=n=1$ on obtient la toute derni\`ere
\'egalit\'e. 
\end{preuve}
Appliquant ce r\'esultat au $\RR$-module $R_e$ on a
\begin{Corollary}\label{cor:modelem}
 L'image de $P$ est contenue dans le centre (additif) de $R_e$. Pour 
$r,s\in R_e$ on a $[r,s]=P((s,r)\cdot H(2))$. Le groupe additif de $R_e$
est nilpotent de classe 2. Et $(-r)s=-(rs)+rPH(s)$. \qed
\end{Corollary}
En particulier pour $r=1$ on a $s+(-1)s=PH(s)$.

\begin{Remark}\label{Anormal}
Si $(M,A)$ est un $\RR$-CP-module, tout sous-groupe de $M$ contenant $A$ est normal, d'apr\`es les relations \ref{pr:modelem}(2) et (MC7b).
\end{Remark}

\noindent\begin{minipage}[t]{0.48\textwidth}
 \begin{Definition}\label{df:BHPsousmod}
Un \emph{$\RR$-BHP-sous-module} de $M$ est un sous-groupe $N$ de $M$ stable
sous les actions de $R_e$ et de $R_{ee}$. Il sera dit \emph{normal} si en
outre $[m,n]\cdot x$ est dans $N$ pour $n\in N$, $m\in M$ et 
$x\in R_{ee}$. 
\end{Definition}

\end{minipage}
\quad 
\begin{minipage}[t]{0.48\textwidth}
 \begin{Definition}\label{df:CPsousmod}
 Un \emph{$\RR$-CP-sous-module} de $(M,A)$ est un $\RR$-CP-module
$(N,B)$  o\`u $N$ (resp. $B$) est un sous-groupe de $M$ (resp. $A$) stable
sous les actions de $R_e$ et de $R_{ee}$. Il sera dit \emph{normal} si en
outre $B=N\cap A$ et si $[m,n]\cdot x$ est dans $N$ pour $n\in N$, $m\in M$
et $x\in R_{ee}$. 
\end{Definition}

\end{minipage}

\smallskip

\noindent\begin{minipage}[t]{0.48\textwidth}
 \begin{Proposition}\label{pr:BHPquotientMC}
 Soient $M$ un $\RR$-BHP-module et $N$ un $\RR$-BHP-sous-module normal de
$M$. Le groupe sous-jacent de $N$ est un sous-groupe normal de $M$, et 
le quotient $M/N$ est naturellement muni d'une structure de
$\RR$-BHP-module.
\end{Proposition}
\end{minipage}
\quad 
\begin{minipage}[t]{0.48\textwidth}
 \begin{Proposition}\label{pr:CPquotientMC}
Soient $\underline{M}=(M,A)$ un $\RR$-CP-module et $\underline{N}=(N,B)$
un $\RR$-CP-sous-module normal de $\underline{M}$. Le groupe sous-jacent de
$N$ est un sous-groupe normal de $M$, et le quotient $(M/N,A/B)$ est
naturellement muni d'une structure de $\RR$-CP-module.
 \end{Proposition}
\end{minipage}
\smallskip
\begin{preuve}[Preuve de \ref{pr:BHPquotientMC}]
 Soit $m\in M$ et $n\in N$, puisque $[n,m]=[m,n]\cdot H(2)\in N$, $N$ est
un sous-groupe normal de $M$. On v\'erifie que les actions de $R_e$ et de
$R_{ee}$ sur $M$ passent bien au quotient:
\begin{equation*}
 (m+n)\cdot r=m\cdot r+n\cdot r+[m,n]\cdot H(r)
\end{equation*}
o\`u les deux derniers termes sont dans $N$.
\begin{equation*}
 [m+n,m']\cdot x=[m,m']\cdot x+[n,m']\cdot x,
\end{equation*}
o\`u le dernier terme est dans $N$.
\end{preuve}
\begin{preuve}[Preuve de \ref{pr:CPquotientMC}] Il reste \`a v\'erifier
les axiomes \ref{ax:MC0}, \ref{ax:MC7a} et \ref{ax:MC7b}.
 Puisque $B=N\cap A$, il est stable sous l'action de $R_e$, et on a
$A/B\simeq A+N/N$ et $(M/N)/(A/B)\simeq M/(A+N)$. Soient alors $m,m'\in M$
et $x\in R_{ee}$. Puisque $[m,m']\cdot x\in A$, sa classe modulo $A+N$ est
nulle, donc $[\ovl m,\ovl{m'}]\cdot x\in A/B$ (o\`u $\ovl{m}$ est la
classe de $m$ modulo $N$). Et si $a\in A$ on a $[m,a]\cdot x=0$. Donc
$[\ovl{m},\ovl{a}]\cdot x=0$.
\end{preuve}

\begin{Proposition}\label{pr:intersecBHP}
 Soit $M$ un $\RR$-BHP-module, et soit $(N_i)_{i\in I}$ une famille de
$\RR$-BHP-sous-modules de $N$. L'intersection $N=\cap_{i\in I}N_i$ est un
$\RR$-BHP-sous-module de $M$.\qed
\end{Proposition}

\begin{Definition}\label{df:BHPssmodeng}
 Soit $M$ un $\RR$-BHP-module et soit $E$ un sous ensemble de $M$ on
appelle $\RR$-BHP-sous-module engendr\'e par $E$ l'intersection des
sous-modules de $M$ qui contiennent $E$. C'est le plus petit sous-module de
$M$ qui contienne $E$.
\end{Definition}
\begin{Proposition}\label{pr:BHPssmodeng}
 Soit $M$ un $\RR$-BHP-module et soit $E$ un sous ensemble de $M$. Le
$\RR$-BHP-sous-module engendr\'e par $E$ est form\'e de l'ensemble des
\'el\'ements de $M$ de la forme 
\begin{equation}
 \sum m_i\cdot r_i +\sum [m'_j,m''_k]\cdot x_{j,k}
\end{equation} 
o\`u $m_i,$ $m'_j$, $m''_k\in E$, $r_i\in R_e$ et $x_{j,k}\in R_{ee}$.\qed
\end{Proposition}

\begin{Definition}\label{df:Rcentre}
 Soit $M$ un $\RR$-BHP-module, un \'el\'ement $m\in M$ est dit
\emph{$\RR$-central} si $[m,n]\cdot x=0$ pour tous $n\in M$ et $x\in
R_{ee}$. 

On appelle \emph{$\RR$-centre} de $M$ l'ensemble $\ZZ_{\RR}(M)$ des $m\in
M$ qui sont $\RR$-centraux.  On appelle \emph{module $\RR$-d\'eriv\'e} de M
le $\RR$-BHP-sous-module $[M,M]_{\RR}$ de $M$, (ou $[M]$ s'il n'y a pas
ambiguit\'e) engendr\'e par les $[m,n]\cdot x$ pour $m,n\in M$ et $x\in
R_{ee}$.
\end{Definition}

\begin{Proposition}\label{pr:Rcentre}
1) Le $\RR$-centre et le module $\RR$-d\'eriv\'e sont des
sous-$\RR$-BHP-modules normaux, et on a les inclusions:
\begin{equation}
 [M,M]\subseteq [M,M]_{\RR}\subseteq \ZZ_{\RR}(M)\subseteq \ZZ(M).
\end{equation}
2) Le $\RR$-centre, le module $\RR$-d\'eriv\'e, le quotient
$M/\ZZ_{\RR}(M)$ et le $\RR$-ab\'elianis\'e 
$M^{\RR\ti ab}=\ovl M=M/[M,M]_{\RR}$ sont des $\ovl R$-modules. De m\^eme, si $(M,A)$ est un $\RR$-CP-module, $A$ et $\ovl{M}=M/A$ sont des $\ovl R$-modules.
\end{Proposition}

\begin{preuve} 
 La premi\`ere assertion d\'ecoule de \ref{ax:MC5} et \ref{ax:MC7}. Les
trois inclusions sont des cons\'equences imm\'ediates de la proposition
\ref{pr:modelem}. On notera en particulier que l'inclusion
$[M,M]_{\RR}\subseteq \ZZ_{\RR}(M)$ \'equivaut \`a la propri\'et\'e
\ref{ax:MC7}.

La seconde vient de ce que $\im P$ annule les quatre modules.
en effet, d'apr\`es (\ref{ax:MC3}) $m\cdot P(x)=[m,m]\cdot x$ qui est dans
le groupe $\RR$-d\'eriv\'e pour $m\in M$, et est nul si $m$ est dans le
$\RR$-centre.
\end{preuve}

\begin{Corollary}\label{cor:BHPvsCPmod}
 Si $M$ est un $\RR$-BHP-module \`a droite, alors $(M,A)$ est un
$\RR$-CP-module \`a droite pour tout sous-groupe $A$ de $M$ stable sous les
actions de $R_e$ et tel que $[M]\subseteq A\subseteq \ZZ_{\RR}(M)$.
Autrement dit, $A$ est un sous-$\bar R$-module de $\ZZ_{\RR}(M)$ contenant
$[M]$.
\end{Corollary}
\begin{preuve} 
 Les propri\'et\'es [\ref{ax:MC7a}] et [\ref{ax:MC7b}] \'equivalent
respectivement \`a $A\subseteq \ZZ_{\RR}(M)$ et \`a $[M,M]_{\RR}\subseteq
A$.
\end{preuve}
On notera que $A$ est stable sous les actions de $R_{ee}$ puisque celles-ci
sont triviales sur $A$.

En particulier, $(R_e,A)$ est un $\RR$-CP-module \`a droite pour tout
sous-groupe $A$ de $R_e$ stable sous les actions de $R_e$ et tel que
$[R_e]\subseteq A\subseteq \ZZ_{\RR}(R_e)$.

\begin{Remark}
 Si $(M,A)$ est un $\RR$-CP-module, il en r\'esulte que 
$[M,M]\subseteq A \subseteq \ZZ(M)$, par cons\'equent le couple $(M,A)$
est un objet de la cat\'egorie $\cp$ d\'efinie en \cite{QmG}.
\end{Remark}


\noindent\begin{minipage}[t]{0.48\textwidth}
 \begin{Proposition}\label{pr:BHPalgsurop}
 Soit $M$ un $\RR$-BHP-module. Alors les groupes ab\'eliens $\mathbb N^*$-gradu\'es
 \[{\rm Gr}_{\gamma}(M) = (M^{\RR\ti ab},[M,M]_{\RR}, 0, \ldots)\]
  \[{\rm Gr}_{ \ZZ}(M) = (M/ \ZZ_{\RR}(M), \ZZ_{\RR}(M), 0, \ldots)\]
   ont une structure naturelle d'alg\`ebre  gradu\'ee  sur l'op\'erade ${\rm OP}(\RR)$.
 \end{Proposition}
\end{minipage}
\quad 
\begin{minipage}[t]{0.48\textwidth}
 \begin{Proposition}\label{pr:CPalgsurop}
Soit $\underline{M}=(M,A)$ un $\RR$-CP-module. Alors le groupe ab\'elien $\mathbb N^*$-gradu\'e
 \[{\rm Gr}(M,A) = (M/A,A, 0, \ldots)\]
a une structure naturelle d'alg\`ebre gradu\'ee sur l'op\'erade  ${\rm OP}(\RR)$.

 \end{Proposition}
\end{minipage}
\smallskip
\begin{preuve}
(pour $(M,A)$, les autres \'etant des cas particuliers) : $M/A$ et $A$ sont des $\ovl{R}$-modules d'apr\`es la proposition \ref{pr:Rcentre}. Vu la graduation et le fait que OP$(\RR)$ est nilpotente de classe $2$, la structure de  OP$(\RR)$-alg\`ebre sur Gr$(M,A)$ est enti\`erement d\'etermin\'ee par les  applications
\[(M/A \oplus A) \otimes {\rm OP}(\RR)_1 \to
M/A \oplus A,\quad (\ovl{m},a)\otimes \ovl{r}\mapsto (\ovl{mr},ar) \]
\[(M/A) \otimes (M/A) \otimes {\rm OP}(\RR)_2 \to A,\quad \ovl{m}\otimes \ovl{m'}\otimes x \mapsto [m,m']\cdot x\,.\]
\end{preuve}

En effet, les propositions pr\'ec\'edentes g\'en\'eralisent le fait que les quotients de la suite centrale descendante d'un groupe (nilpotent de classe $2$ ici) forment une alg\`ebre de Lie, cf.\ le paragraphe \ref{sec:Pnul} plus loin.


\subsection{Morphismes de modules sur un anneau carr\'e}
\label{sec:morph}

Soit $M$ et $N$ deux $\RR$-BHP-modules, une application $f$ de $M$ dans $N$
sera dite additive (et non pas lin\'eaire) si elle respecte l'addition,
elle sera dite $R_e$-\'equivariante si elle respecte l'action de $R_e$, et
elle sera dite $R_{ee}$-\'equivariante si elle respecte l'action de
$R_{ee}$ sur $M$ et $N$, c'est \`a dire si
\begin{equation*}
 f([m,m']\cdot x)=[f(m),f(m')]\cdot x.
\end{equation*}
Elle sera dite $R_e$-lin\'eaire si elle est additive et
$R_e$-\'equivariante.
\smallbreak
\noindent\begin{minipage}[t]{0.48\textwidth}
 \begin{Definition}(\cite{BHP})\label{df:morBHPmodc}
 Soient $M$ et $N$ deux $\RR$-BHP-modules, un BHP-mor\-phisme ou
application $\RR$-BHP-lin\'eaire est une application $\alpha:M\to N$ qui
est $R_e$-lin\'eaire et $R_{ee}$-\'equivariante.
\end{Definition}
\end{minipage}
\quad 
\begin{minipage}[t]{0.48\textwidth}
 \begin{Definition}(\cite{chap2})\label{df:morCPmodc}
 Soient $(M,A)$ et $(N,B)$ deux $\RR$-CP-modules un $\RR$-CP-morphisme de
$(M,A)$ dans $(N,B)$ est
une application $\alpha:M\to N$ qui est $R_e$-lin\'eaire et
$R_{ee}$-\'equivariante, et telle que $\alpha(A)\subseteq B$.
\end{Definition}
\end{minipage}

\smallbreak
L'application $\alpha$  est donc un homomorphisme de groupe tel que 
$\alpha(m\cdot r)=\alpha(m)\cdot r$, et 
$\alpha([m,m']\cdot x)=[\alpha(m),\alpha(m')]\cdot x$.

\smallbreak
\noindent\begin{minipage}[t]{0.48\textwidth}
 Il est clair qu'on obtient ainsi une cat\'egorie not\'ee BHP-$\modr$. Si
$N$ est un sous $\RR$-BHP-module normal de $M$ on peut v\'erifier que $N$
est le noyau dans BHP-$\modr$ du morphisme $M\to M/N$.
\end{minipage}
\quad 
\begin{minipage}[t]{0.48\textwidth}
 Il est clair qu'on obtient ainsi une cat\'egorie not\'ee CP-$\modr$. Si
$(N,B)$ est un sous $\RR$-CP-module normal de $(M,A)$ on peut v\'erifier
que $(N,B)$ est le noyau dans CP-$\modr$ du morphisme 
$(M,A)\to (M/N,A/B)$.
\end{minipage}
\smallbreak

Il est imm\'ediat qu'\`a un $\RR$-BHP-morphisme $M\to N$ on associe
canoniquement un $\RR$-CP-morphisme $(M,[M])\to (N,[N])$, ou plus
g\'en\'eralement $(M,[M])\to (N,B)$ o\`u $B$ est un sous-BHP-module tel que
$[N]\subseteq B\subseteq \ZZ_{\RR}N$.

\begin{Example}\label{ex:mulgauche}
 Soit $M$ un $R$-BHP module, pour tout $m\in M$  
l'application $g_m: R_e\to M$, $r\mapsto m\cdot r$ est
$\RR$-BHP-lin\'eaire. En effet, elle est additive et $R_e$ \'equivariante
par \ref{ax:MC1}, et $R_{ee}$-\'equivariante puisque, d'apr\`es
\ref{ax:MC3} et \ref{ax:MC6} on a:
\begin{equation*}
 m\cdot([r,s]\cdot x)=m\cdot P((r,s)\cdot x)=[m,m]\cdot((r,s)\cdot x)
=[m\cdot r,m\cdot s]\cdot x
.\end{equation*}
En particulier, si $M=R_e$ les multiplications \`a gauche dans $R_e$ sont
$\RR$-BHP-lin\'eaires.
 \end{Example}


\smallbreak
\noindent\begin{minipage}[t]{0.48\textwidth}
 \begin{Proposition}\label{pr:GrmapBHP}
Une application $\RR$-BHP-lin\'eaire $f\colon M \to N$ entre $\RR$-BHP-modules $M,N$ induit un morphisme  de ${\rm OP}(\RR)$-alg\`ebres gradu\'ees
\[ {\rm Gr}_{\gamma}(f) \colon {\rm Gr}_{\gamma}(M) \to {\rm Gr}_{\gamma}(N)\, \]
de fa\c{c}on \'evidente.
 \end{Proposition}
\end{minipage}
\quad 
\begin{minipage}[t]{0.48\textwidth}
 \begin{Proposition}\label{pr:GrmapCP}
Une application $\RR$-CP-lin\'eaire $f\colon (M,A) \to (N,B) $ entre $\RR$-CP-modules $(M,A),(N,B)$ induit des morphismes de ${\rm OP}(\RR)$-alg\`ebres gradu\'ees
\[ {\rm Gr}(f) \colon {\rm Gr} (M,A) \to {\rm Gr} (N,B) \,,\]
de fa\c{c}on \'evidente.
 \end{Proposition}

\end{minipage}

On obtient ainsi des foncteurs
\[{\rm Gr}_{\gamma} \colon \mbox{BHP-$\modr$} \to {\rm GrAlg}_2({\rm OP}(\RR))\,, \]
\[{\rm Gr} \colon \mbox{CP-$\modr$} \to {\rm GrAlg}_2({\rm OP}(\RR)) \]
o\`u ${\rm GrAlg}_2({\rm OP}(\RR))$ d\'esigne la cat\'egorie des alg\`ebres gradu\'es nilpotentes de classe $2$ sur l'op\'erade ${\rm OP}(\RR)$. Ils g\'en\'eralisent le foncteur qui associe \`a un groupe nilpotent de classe $2$ l'alg\`ebre de Lie form\'ee des quotients de la suite centrale descendante, cf.\ le paragraphe \ref{sec:Pnul} plus loin.


\begin{Definition}\label{df:modCrochZero}
 Un $\RR$-module est dit \emph{$\RR$-ab\'elien}, si l'application $M\times
M\times R_{ee}\to M$ est nulle.
\end{Definition}

\begin{Lemma}\label{lm:modCrochZero}
 Soit $M$ un $\RR$-module $\RR$-ab\'elien. Alors
\begin{MYitemize}
 \item le groupe additif de $M$ est ab\'elien,
\item $M$ est annul\'e par $\im P$,
\item l'action de $R_e$ sur $M$ est distributive \`a gauche, 
\end{MYitemize}
et $M$ est naturellement un $\bar R$-module.
\end{Lemma}
\begin{preuve}
 Cela d\'ecoule imm\'ediatement de \ref{pr:modelem} (2) et des axiomes
(\ref{ax:MC3})  et (\ref{ax:MC2}).
\end{preuve}

\smallbreak
\noindent\begin{minipage}[t]{0.48\textwidth}
 \begin{Corollary}\label{cor:mczBHP}
 La sous cat\'egorie pleine de BHP-$\modr$ form\'ee des modules
$\RR$-ab\'eliens est \'equivalente \`a la cat\'egorie des $\bar R$-modules.
 \end{Corollary}
\end{minipage}
\quad 
\begin{minipage}[t]{0.48\textwidth}
 \begin{Corollary}\label{cor:mczCP}
 La sous-cat\'egorie pleine de CP-$\modr$ form\'ee des modules
$\RR$-ab\'eliens est \'equivalente \`a la cat\'egorie des couples de $\bar
R$-modules.
 \end{Corollary}

\end{minipage}

\smallbreak
\noindent\begin{minipage}[t]{0.48\textwidth}
 \begin{Proposition}\label{pr:factabBHP}
 Soient $M$ et $N$ deux $\RR$-BHP-modules, et supposons que $N$ est
$\RR$-ab\'elien. Toute application $\RR$-BHP-lin\'eaire de $M$ dans $N$ se
factorise par $\bar M=M/[M]_{\RR}$. On obtient donc une bijection naturelle
$$\bar R\ti mod(\bar M,N)\to \text{BHP-}\modr(M,N).\qedhere $$
 \end{Proposition}
\end{minipage}
\quad 
\begin{minipage}[t]{0.48\textwidth}
 \begin{Proposition}\label{pr:factabCP}
 Soient $(M,A)$ et $(N,B)$ deux $\RR$-CP-modules, et supposons que $(N,B)$
est $\RR$-ab\'elien. Toute application $\RR$-CP-lin\'eaire de $(M,A)$ dans
$(N,B)$ se factorise par $(\bar M,\bar A)=(M/[M]_{\RR},A/[M]_{\RR}$. On
obtient donc une bijection naturelle 
\begin{multline*} \bar R\ti mod((\bar M,\bar A), (N,B))\to\\
\text{CP-}\modr((M,A),(N,B)).\qedhere 
\end{multline*}
 \end{Proposition}
\end{minipage}

\smallbreak
\noindent\begin{minipage}[t]{0.48\textwidth}
 \begin{Corollary}\label{cor:morphmczBHP}
  Sous les m\^emes hypoth\`eses, l'ensemble BHP-$\modr(M,N)$ est
naturellement un groupe ab\'elien. Si en outre l'anneau $\bar R$ est
commutatif, c'est un $\bar R$-module. \qed
 \end{Corollary}
\end{minipage}
\quad 
\begin{minipage}[t]{0.48\textwidth}
 \begin{Corollary}\label{cor:morphmczCP}
  Sous les m\^emes hypoth\`eses, l'ensemble CP-$\modr((M,A),(N,B))$ est
naturellement un groupe ab\'elien. Si en outre l'anneau $\bar R$ est
commutatif, c'est un $\bar R$-module. \qed
 \end{Corollary}
\end{minipage}
\smallbreak
\noindent\begin{minipage}[t]{0.48\textwidth}
 \begin{Example}\label{ex:lib1genBHP}
 Soit $M$ un $\RR$-BHP-module et $m\in M$. Il existe une unique
application $\RR$-lin\'eaire $R_e\to M$ qui envoie 1 sur $m$. En effet, si
une telle application $\varphi$ existe, comme elle est
$R_e$-\'equivariante, on a $\varphi(r)=\varphi(1\cdot r) =\varphi(1)\cdot
r=m\cdot r$. d'o\`u l'unicit\'e. On a vu \ref{ex:mulgauche} que $r\mapsto
m\cdot r$ est $\RR$-BHP-lin\'eaire. On en d\'eduit alors que $R_e$ est le
$\RR$-BHP-module libre \`a un g\'en\'erateur.
\end{Example}
\end{minipage}
\quad 
\begin{minipage}[t]{0.48\textwidth}
 \begin{Example}\label{ex:lib1genCP}
 Soit $(M,A)$ un $\RR$-CP-module et $m\in M$. Il existe une
unique application $\RR$-CP-lin\'eaire $(R_e,P(R_{ee}))\to (M,A)$ qui
envoie 1 sur $m$. On en d\'eduit que $(R_e,P(R_{ee}))$ est le
$\RR$-CP-module libre \`a un g\'en\'erateur.
\end{Example}
\end{minipage}
\smallskip

\begin{Definition}\label{df:Rimage}
 Soit $E$ un ensemble, $M$ un $\RR$-BHP-module, et $f:E\to M$ une
application. On appelle $\RR$-image de $f$, qu'on note $\im_{\RR} f$ le
$\RR$-BHP-sous-module de $M$ engendr\'e par l'ensemble $f(E)$ (cf.
\ref{df:BHPssmodeng}).
\end{Definition}

\begin{Lemma}\label{lm:centrImage}
 Avec les notations ci-dessus, si $Y$ est une partie de 
$f(E)$ qui $\RR$-centralise tous les \'el\'ements de $f(E)$, alors le
$\RR$-BHP-sous-module de $M$ engendr\'e par $Y$ est $\RR$-central dans
$\im_{\RR}f$. 
\end{Lemma}
\begin{preuve}
Cela d\'ecoule imm\'ediatement de la proposition \ref{pr:BHPssmodeng}. 
\end{preuve}

\subsection{Exemples}
\label{sec:exBHP}

Nous rappelons ici des exemples du paragraphe 8 de \cite{BHP}.

\subsubsection{Modules sur un anneau commutatif classique}
\label{sec:classique}

Soit $R$ un anneau commutatif classique et $\RR=(R\to 0\to R)$ l'anneau
carr\'e associ\'e. Un BHP-module $M$ est simplement un $R$-module classique
(les propri\'et\'es (\ref{ax:MC3}) \`a (\ref{ax:MC7}) sont vides, et le
groupe $M$ est commutatif puisque $H=0$).

\subsubsection{Modules sur l'anneau $\Znil$}
\label{sec:Znil}
On a $\Znil=\annc{\Z}$ avec $\Z_e=\Z_{ee}=\Z$, $H(n)=\binom{n}{2}$ et
$P=0$.
Un $\Znil$-BHP-module est un groupe nilpotent de classe 2. L'action de
$\Z_e$ est donn\'ee par $m\cdot r=m+\dots m$ ($r$ fois); celle de $\Z_{ee}$
est d\'etermin\'ee par $[m,m']\cdot 1=[m',m]$. 

\subsubsection{Modules sur les anneaux carr\'es tels que $P=0$}
\label{sec:Pnul}
Soit $R$ un anneau, prenons $R_e=R$ et soit $R_{ee}$ un 
$R_e\tp R_e\tp R_e^{op}$-module tel que pour $r,\ s,\ t\in R_e$
\begin{equation}\label{eq:symg}
 (r\tp s)\cdot x\cdot t= (s\tp r)\cdot x \cdot t.
\end{equation}
Pour toute application $H:R_e\to R_{ee}$ v\'erifiant
\begin{align*}
 H(r+s)&=H(r)+H(s)+(s,r)\cdot H(2)\qquad\text{et}&
H(r\,s)&=(r\tp r)\cdot H(s)+H(r)\cdot s
\end{align*}
on obtient un anneau carr\'e
$R_{\nil}^H=(\xymatrix{R_e\ar[r]^H&R_{ee}\ar[r]^{P=0}&R_e})$, o\`u
$T=-\Id$. 
Et tout anneau carr\'e tel que $P=0$ est obtenu de cette fa\c con.

Un $R_{\nil}^H$-module $M$ est alors un groupe muni d'une action \`a droite
du mono\"ide $R_e$, et d'op\'erations index\'ees par $R_{ee}$ v\'erifiant
les axiomes de la d\'efinition \ref{df:BHPmodc}. Les axi\^omes  \ref{ax:MC3} et \ref{ax:MC4}
  entra\^inent  que ces  op\'erations sont toutes altern\'ees.

En particulier, si $R$ est commutatif et si on prend pour $R_{ee}$
l'id\'eal $I_2$ de $R$ engendr\'e par les \'el\'ements $r^2-r$ avec
$H(r)=r^2-r$, on obtient l'anneau carr\'e  $R_{\Nil}=(R\to I_2\to R)$
\cite{GH}.

Lorsque $R$ est un anneau 2-bin\^omial, l'anneau carr\'e 
$R_\nil=(\xymatrix{R\ar[r]^H&R\ar[r]^{P=0}&R})$ avec $H(r):=\binom{r}{2}$
est isomorphe \`a $R_\Nil$. Les $R$-groupes nilpotents de classe $\leq 2$ classiques
(\cite{Wf}) sont donc exactement les $R_{\Nil}$-modules. La notion de
$R_{\Nil}$-module pour un anneau commutatif quelconque g\'en\'eralise
donc la notion de $R$-groupe nilpotent de classe $\leq 2$ pour un anneau
2-bin\^omial. 

Notons aussi que Op$(R_\nil)$ est l'op\'erade Lie tronqu\'ee, et que pour un $R_\nil$-BHP-module $M$ (c\`ad, un $R$-groupe nilpotent de classe $\leq 2$), la Op$(R_\nil)$-alg\`ebre Gr$_{\gamma}(M)$ est la $R$-alg\`ebre de Lie gradu\'ee associ\'ee  \`a $M$
classique.

\subsubsection{Modules sur l'anneau carr\'e $\La_R$ de \cite{BHP} (8.8)}
\label{sec:anncarLie}
Soit $R$ un anneau commutatif, soit
\begin{equation}
 \La_R:=(\xymatrix{R\ar[r]^0&R\ar[r]^0&R})
\end{equation} 
l'anneau carr\'e o\`u $R_e=R$ avec les m\^emes addition et multiplication,
$R_{ee}=R$ avec la m\^eme addition, $H=P=0$. Les actions de $R_e$ sur
$R_{ee}$ sont donn\'ees par :
$(r,s)\cdot x\cdot t:= rsxt$.

Soit maintenant $M$ un $\La_R$-BHP-module. D'apr\`es \ref{pr:modelem} (2)
le groupe additif de $M$ est commutatif, et d'apr\`es \ref{ax:MC2} $M$ est
un $R$-module (au sens classique). D'apr\`es \ref{ax:MC6} les crochets sur
$M$ sont d\'etermin\'es par le seul crochet $[m,n]\cdot 1$ qui est
bilin\'eaire en $m$ et $n$ d'apr\`es \ref{ax:MC5} et qui doit v\'erifier
$[m,m]\cdot 1=0$, $[n,m]\cdot 1=-[m,n]\cdot 1$ et $[[m,m']\cdot 1,n]\cdot
1]=0$. Il en r\'esulte qu'un $\La_R$-BHP-module n'est rien d'autre qu'une
$R$-alg\`ebre de Lie nilpotente de classe $\leq 2$.

On v\'erifie de m\^eme qu'un $\La_R$-CP-module est un couple $(M,A)$, o\`u
$M$ est une $R$-alg\`ebre de Lie nilpotente de classe $\leq 2$, et $A$ un
sous module de $M$ contenant $[M,M]$ et tel que $[A,M]=0$.

\subsubsection{Modules sur l'anneau carr\'e $\bigotimes_R$ de \cite{BHP}
(8.8)}
\label{sec:anncarTenseur}
Dans cet exemple et les deux qui vont suivre, soit $R$ un anneau
commutatif, et soit $R_e$ le groupe $R\oplus R$ muni de la multiplication
\begin{equation}\label{eq:multRe}
 (r,s)(r',s'):=(rr',r^2s'+sr').
\end{equation} 
On v\'erifie que $R_e$ satisfait tous les axiomes d'un anneau, sauf la
distributivit\'e \`a droite. On a en effet:
\begin{align*}
 (r+r',s+s')(r'',s")&=(rr''+r'r'',(r+r')^2s''+(s+s')r''),\\
(r,s)(r'',s'')+(r',s')(r'',s'') &=(rr'',r^2s''+sr'')+
(r'r'',r'^2s''+s'r''),\\
&=(rr''+r'r'',r^2s''+sr''+r'^2s''+s'r'').
\end{align*}
Donc
\begin{equation*}
 ((r,s)+(r',s'))(r'',s'')=(r,s)(r'',s'')+ (r',s')(r'',s'')+(0,2rr's'').
\end{equation*}
On remarquera que l'inclusion canonique $R\to R_e$, $r\mapsto (r,0)$ est
compatible avec les deux op\'erations.

D\'efinissons alors un anneau carr\'e:
\begin{equation}\label{def:tensR}
\textstyle{\bigotimes_R}:=(\xymatrix{R_e\ar[r]^{H\quad }& R\oplus
R\ar[r]^{\quad P}&R_e})
\end{equation} 
o\`u $H(r,s):=(s,s)$, $P(x,y):=(0,x+y)$ et o\`u les actions de $R_e$ sur
$R_{ee}$ sont donn\'ees par
\begin{equation*}
 ((r,s),(r',s'))\cdot(x,y)\cdot(r'',s''):= (rr'xr'',rr'yr'').
\end{equation*}
Le lecteur v\'erifiera alors que $\bigotimes_R$ satisfait les axiomes d'un
anneau carr\'e. Par exemple, pour \ref{ax:AC7} on a 
\begin{align*}
 P(((r,s),(r'',s'))\cdot H(r'',s''))&=
 P(((r,s),(r',s'))\cdot (s'',s'')=
P((rr's'',rr's''))=(0,2rr's'').
\end{align*}
On a $T(x,y)=HP(x,y)-(x,y)=H(0,x+y)-(x,y)= (x+y,x+y)-(x,y)=(y,x)$.

Soit maintenant $M$ un $\otimes_R$-BHP-module. Comme on a
$H(2)=H((2,0))=0$, la propri\'et\'e \ref{pr:modelem} (2) montre que le
groupe additif de $M$ est commutatif. En outre, $M$ est naturellement un
$R$-module via l'inclusion de $R$ dans $R_e$.

D'apr\`es \ref{ax:MC5} \ref{ax:MC4} et \ref{ax:MC6} on a 
\begin{align*}
 [m,n]\cdot(r,s)&=[m,n]\cdot(r,0)+[m,n]\cdot(0,s)=
[m,n]\cdot(r,0)+[m,n]\cdot T(s,0)=[m,n]\cdot(r,0)+[n,m]\cdot(s,0)\\
&= ([m,n]\cdot (1,0))r+([n,m]\cdot (1,0))s.
\end{align*}
Les crochets sont donc enti\`erement d\'etermin\'es par le seul crochet
$[m,n]\cdot 1$.  Si on pose $m*n:=[m,n]\cdot 1$, le calcul pr\'ec\'edent
s'\'ecrit:
\begin{equation*}
 [m,n]\cdot(r,s)=(m*n)r+(n*m)s.
\end{equation*}
L'action de la seconde composante de $R_e$ sur $M$ s'exprime aussi \`a
l'aide de l'op\'eration $*$ gr\^ace \`a \ref{ax:MC3}:
\begin{equation*}
 m\cdot (0,s)=m\cdot P(s,0)=[m,m]\cdot (s,0)=(m*m)s
\end{equation*}
Il en r\'esulte que la donn\'ee d'un $\bigotimes_R$-BHP-module $M$
\'equivaut \`a la donn\'ee d'un $R$-module $M$ muni d'une multiplication
$*$ qui est $R$-bilin\'eaire et nilpotente de classe $\leq 2$, puisque
\ref{ax:MC7} avec $x=y=1$ donne $(m*m')*n=0$ pour tous $m$, $m'$ et $n$.
Autrement dit, $M$ est une $R$-alg\`ebre non n\'ecessairement commutative
(sans unit\'e) nilpotente de classe $\leq2$.

On v\'erifie de m\^eme qu'un $\bigotimes_R$-CP-module est un couple
$(M,A)$, o\`u $M$ est une $R$-alg\`ebre non n\'ecessairement commutative
(sans unit\'e) nilpotente de classe $\leq 2$, et $A$ un id\'eal de $M$
contenant $M*M$ et
tel que $A*M=M*A=0$.

\subsubsection{Modules sur l'anneau carr\'e $S_R$ de \cite{BHP} (8.8)}
\label{sec:anncarS}
A partir du m\^eme $R_e$, d\'efinissons un autre anneau carr\'e:
\begin{equation}\label{def:symR}
\textstyle{S_R}:=(\xymatrix{R_e\ar[r]^{H}& R\ar[r]^{P}&R_e})
\end{equation} 
o\`u $H(r,s):=2s$, $P(x):=(0,x)$ et o\`u les actions de $R_e$ sur $R_{ee}$
sont donn\'ees par
\begin{equation*}
 ((r,s),(r',s'))\cdot x\cdot(r'',s''):= rr'xr''.
\end{equation*}
Le lecteur v\'erifiera alors que $S_R$ satisfait les axiomes d'un anneau
carr\'e. Par exemple, pour \ref{ax:AC7} on a 
\begin{align*}
 P(((r,s),(r',s'))\cdot H(r'',s''))&=
 P(((r,s),(r',s'))\cdot (2s''))=
P((2rr's''))=(0,2rr's").
\end{align*}
On a $T(x)=HP(x)-x=H(0,x)-x= 2x-x=x$.

Soit maintenant $M$ un $S_R$-BHP-module. Comme on a $H(2)=H((2,0))=0$, la
propri\'et\'e \ref{pr:modelem} (2) montre que le groupe additif de $M$ est
commutatif. En outre, $M$ est naturellement un $R$-module via l'inclusion
de $R$ dans $R_e$.

D'apr\`es \ref{ax:MC6} on a 
\begin{equation*}
 [m,n]\cdot x=[m,n]\cdot(1\cdot (x,0))= ([m,n]\cdot 1)(x,0)=([m,n]\cdot
1)x.
\end{equation*}
Les crochets sont donc enti\`erement d\'etermin\'es par le seul crochet
$[m,n]\cdot 1$. D\'efinissons alors une multiplication sur $M$ par
$m*n:=[m,n]\cdot 1$. 

L'action de la seconde composante de $R_e$ sur $M$ s'exprime aussi \`a
l'aide de l'op\'eration $*$ gr\^ace \`a \ref{ax:MC3}:
\begin{equation*}
 m\cdot (0,s)=m\cdot P(s)=[m,m]\cdot s=(m*m)s
\end{equation*}

En outre, d'apr\`es \ref{ax:MC4} on a $[m,n]\cdot T(1)=[n,m]\cdot 1$. Et
comme $T=\Id$ on en d\'eduit que la multiplication * est commutative.

Il en r\'esulte que la donn\'ee d'un $S_R$-BHP-module $M$ \'equivaut \`a la
donn\'ee d'un $R$-module $M$ muni d'une multiplication commutative $*$ qui
est $R$-bilin\'eaire et nilpotente de classe $\leq 2$, puisque \ref{ax:MC7}
avec $x=y=1$ donne $(m*m')*n=0$ pour tous $m$, $m'$ et $n$. Autrement dit,
$M$ est une $R$-alg\`ebre commutative (sans unit\'e) nilpotente de classe
$\leq2$.

On v\'erifie de m\^eme qu'un $S_R$-CP-module est un couple $(M,A)$, o\`u
$M$ est une $R$-alg\`ebre commutative nilpotente de classe $\leq 2$, et $A$
un id\'eal de $M$ contenant $M*M$ et tel que $A*M=0$.

\subsubsection{Modules sur l'anneau carr\'e $\Ga_R$ de \cite{BHP} (8.8)}
\label{sec:anncarGamma}
Toujours \`a partir du m\^eme $R_e$ d\'efinissons un troisi\`eme anneau
carr\'e:
\begin{equation}\label{def:GammaR}
\textstyle{\Ga_R}:=(\xymatrix{R_e\ar[r]^{H}& R\ar[r]^{P}&R_e})
\end{equation} 
o\`u $H(r,s):=s$, $P(x):=(0,2x)$ et o\`u les actions de $R_e$ sur $R_{ee}$
sont donn\'ees par
\begin{equation*}
 ((r,s),(r',s'))\cdot x\cdot(r'',s''):= rr'xr''.
\end{equation*}
Le lecteur v\'erifiera alors que $\Ga_R$ satisfait les axiomes d'un anneau
carr\'e. Par exemple, pour \ref{ax:AC7} on a 
\begin{align*}
 P(((r,s),(r',s'))\cdot H(r'',s''))&=
 P(((r,s),(r',s'))\cdot (s''))=
P((rr's''))=(0,2rr's").
\end{align*}
On a $T(x)=HP(x)-x=H(0,2x)-x= 2x-x=x$.

Soit maintenant $M$ un $\Ga_R$-BHP-module. Comme on a $H(2)=H((2,0))=0$, la
propri\'et\'e \ref{pr:modelem} (2) montre que le groupe additif de $M$ est
commutatif. En outre, $M$ est naturellement un $R$-module via l'inclusion
de $R$ dans $R_e$. D\'efinissons alors une application $\ga:M\to M$ par
$\ga(m):=m\cdot(0,1)$. D'apr\`es \ref{ax:MC2} on a
\begin{equation*}
 [m,n]\cdot 1=[m,n]\cdot H(0,1)=(m+n)\cdot(0,1)-m\cdot(0,1)-n\cdot(0,1)=
\ga(m+n)-\ga(m)-\ga(n).
\end{equation*}
Il en r\'esulte que la multiplication $m*n:=[m,n]\cdot(0,1)$ est
commutative et $R$-bilin\'eaire, et que $m*m=2\ga(m)$. On a aussi
\begin{equation*}
 \ga(mr)=\ga(m\cdot(r,0))=m\cdot(r,0)\cdot(0,1)=m\cdot(0,r^2)=m\cdot(0,
1)\cdot (r^2,0)=\ga(m)\cdot(r^2,0)=\ga(m)r^2.
\end{equation*}
Il en r\'esulte que
la donn\'ee d'un $\Ga_R$-BHP-module $M$ \'equivaut \`a la donn\'ee d'une
$R$-alg\`ebre \`a puissance divis\'ee (sans unit\'e) nilpotente de
classe $\leq 2$. 

\subsubsection{Anneaux carr\'es associ\'es \`a la suspension des
pseudo-plans projectifs \cite{BHP} (8.4)}
\label{sec:anncarSigPn}
Dans l'exemple pr\'ec\'edent, soit $I_2$ l'id\'eal de $R$ engendr\'e par
les \'el\'ements $r^2-r$ et soit $\varepsilon\in Ann(I_2)$. D\'eformons
l'addition de $R_e$ en posant
\begin{equation*}
 (r,s)+(r',s'):=(r+r',s+s'+\varepsilon\, r\,r')
\end{equation*}
Avec le m\^eme $R_{ee}$, les m\^emes actions de $R_e$ sur $R_{ee}$ et les
m\^emes applications $H$ et $P$ que pour l'anneau carr\'e $\Ga_R$, on
obtient un nouvel anneau carr\'e $\Ga_R^{\varepsilon}$, qui g\'en\'eralise
l'exemple pr\'ec\'edent puisque $\Ga_R^0=\Ga_R$.

Lorsque l'on prend $R=\Z/n\Z$ et $\varepsilon=\binom{n}{2}$ l'anneau
carr\'e obtenu est l'anneau carr\'e des endomorphismes de $\Sigma P_n$.

\section{Anneaux carr\'es commutatifs et applications $\RR$-quadratiques}
\label{sec:ACCM}

\subsection{Anneaux carr\'es commutatifs}
\label{sec:AnnCarComm}

\begin{Definition}\label{df:accomm}
Un anneau carr\'e $\RR$ est dit commutatif si
\begin{MYitemize}
 \item l'anneau quotient $\ovl R=R_e/\im P$ est commutatif
 \item les trois actions de $R_e$ sur $R_{ee}$ co\"incident,(\footnote{Une
autre d\'efinition possible a finalement \'et\'e \'ecart\'ee : Un anneau
carr\'e $\RR$ est dit commutatif si la multiplication dans $R_e$ est
commutative et si les trois actions de $R_e$ sur $R_{ee}$  co\"incident. En
effet la commutativit\'e de la multiplication de $R_e$ entra\^ine la
distributivit\'e \`a gauche, il en r\'esulte que $PH=0$ et que $2P=0$. Elle
entra\^ine aussi la commutativit\'e de l'addition de $R_e$, puisque
$[m,n]=P((n,m)\cdot H(2))=PH(2)nm=0$. On aurait finalement que $R_e$ serait
un anneau commutatif, ce qui semble un peu trop restrictif.}) c'est \`a
dire $(r,1)\cdot x=(1,r)\cdot x=x\cdot r$.
\end{MYitemize}
\end{Definition}


Remarquons que la propri\'et\'e de commutativit\'e s'exprime en termes de l'op\'erade associ\'ee : $\RR$ est commutatif si et seulement si l'anneau ${\rm OP}(\RR)_1$ est commutatif et 
${\rm OP}(\RR)$ est une op\'erade de ${\rm OP}(\RR)_1$-modules.


\begin{Lemma}\label{lm:quasicomm}
 Si $\RR$ est commutatif, on a $x\cdot rs=x\cdot sr$ pour $r,s\in R_e$ et 
$x\in R_{ee}$. En outre le groupe additif de $R_e$ est commutatif.
\end{Lemma}
\begin{preuve}
La premi\`ere relation est une cons\'equence imm\'ediate de la
d\'efi\-nition.
Pour la seconde propri\'et\'e, soient $r$ et $s$ dans $R_e$, on a d'apr\`es
\ref{pr:modelem}
\begin{equation*}
 [r,s]=P((s,r)\cdot H(2))=P((1,1)\cdot H(2))\,rs=[1,1]\,rs=0.\qedhere
\end{equation*}

\end{preuve}

\begin{Lemma}\label{lm:basrelmod}
 Soit $M$ un module sur l'anneau carr\'e commutatif $\RR$. Pour $m,n\in M$,
$r,s,t\in R_e$ et $x,y\in R_{ee}$, on a
$ [m\cdot r,n\cdot s]\cdot (x\cdot t)=([m,n]\cdot x)\cdot rst$,
et dans le produit $rst$ on peut changer l'ordre des facteurs.
\end{Lemma}

\begin{preuve}
 Cela d\'ecoule de  \ref{ax:MC6}.
\end{preuve}

\begin{Corollary}\label{cor:Rmult}
 Soit $\RR$ un anneau carr\'e commutatif, pour $r\in R_e$ et $x\in R_{ee}$
on a $rP(x)=P(x)r^2$.
\end{Corollary}
\begin{preuve}
 En utilisant \ref{ax:AC1}, la d\'efinition \ref{df:accomm} et \`a nouveau
\ref{ax:AC1} on a:
\begin{equation*}
 rP(x)=P((r,r)\cdot x)=P(x\cdot r^2)=P(x)r^2.
\end{equation*}

\end{preuve}

\begin{Remark}\label{rm:pseudocom}
 Puisque $\ovl R$ est commutatif, $rs-sr\in \im P$, par cons\'equent 
$m\cdot sr-m\cdot rs= m\cdot(sr-rs)=m\cdot P(x)=[m,m]\cdot x$. Donc, si $m$
est tel que $[m,m]\cdot x=0$ pour tout $x\in R_{ee}$, on a $m\cdot
sr=m\cdot rs$ pour tous $r,s\in R_e$. 
\end{Remark}

\subsubsection{Exemples}
\label{sec:exemples}

Les anneaux carr\'es des exemples
\ref{sec:classique}, $\Z_{\nil}$ (\ref{sec:Znil}), $R_{\nil}^H$ si $R$ est
commutatif (\ref{sec:Pnul}) et si $(r,1)\cdot x=x\cdot r$,
$\La_R$ (\ref{sec:anncarLie}), $\bigotimes_R$ (\ref{sec:anncarTenseur}) et
$S_R$ (\ref{sec:anncarS}) sont des anneaux carr\'es commutatifs.

\subsubsection{Les anneaux carr\'es $\Ga_R$ (\ref{sec:anncarGamma}) et
$\Ga_R^{\varepsilon}$ (\ref{sec:anncarSigPn})}
Pour ces anneaux carr\'es, on a $\im P=2R$, donc $\bar R=R\times R/2R$. On
v\'erifie alors que $\bar R$ est un  anneau commutatif si et seulement si
$I_2=2R$. Cette condition est r\'ealis\'ee, par exemple si $1/2\in R$ ou
si $R$ est un anneau 2-bin\^omial, ou si $R=\Z/n\Z$.


\subsection{Applications quadratiques entre $\RR$-modules}
\label{sec:appRq}

\textit{Dans toute la suite, on supposera que $\RR$ est un anneau carr\'e commutatif.}

\begin{Definition}\label{df:defauts}
 Soient $\RR$ un anneau carr\'e commutatif, $M$ et $N$ deux $\RR$-modules
et $f:M\to N$ une application. On appelle \emph{d\'efaut additif} de $f$
l'application
\begin{align}
 d_f:M\times M&\to N&(m,m')\mapsto d_f(m,m')&:=f(m+m')-f(m')-f(m).
\end{align} 
On appelle \emph{d\'efauts scalaires} de $f$ les applications
\begin{align}
 f_{(r)}:M&\to N& m\mapsto f_{(r)}(m)&:=f(m\cdot r)-f(m)\cdot r,
\end{align}
pour $r\in R_e$. On appelle \emph{d\'efauts des commutateurs} de $f$ les
applications
\begin{align}
f_{[x]}:M\times M&\to N& (m,m')\mapsto f_{[x]}(m,m')&:=f([m,m']\cdot x)
-[f(m),f(m')]\cdot x,
\end{align}
pour $x\in R_{ee}$.
\end{Definition}

\noindent\begin{minipage}[t]{0.48\textwidth}
 \begin{Definition}\label{df:appBHPqu}
 Soit $\RR$ un anneau carr\'e commutatif et soient $M$ et $N$ deux
$\RR$-BHP-modules, une application $f:M\to N$ est
$\RR$-BHP-\emph{quadratique} si elle satisfait:
\begin{enumerate}
 \item les images des d\'efauts : $\im_{\RR} d_f$, $\im_{\RR} f_{(r)}$ et
$\im_{\RR} f_{[x]}$ sont $\RR$-centrales dans $\im_{\RR} f$, pour $r\in
R_e$ et $x\in R_{ee}$.
 \item le d\'efaut additif $d_f$ et les d\'efauts des commutateurs 
$f_{[x]}$  sont $\RR$-bilin\'eaires,
\item les d\'efauts scalaires $f_{(r)}$ sont homog\`enes de degr\'e 2
(\footnote{\ Nous n'utilisons pas ici l'appellation "quadratique
homog\`ene" car pour le moment rien ne dit que $f_{(r)}$ est
$\RR$-quadratique. On verra plus loin \ref{pr:defscalquadBHP} que si $f$
est $\RR$-BHP-quadratique, les $f_{(r)}$ le sont aussi.}), i.e.
$f_{(r)}(m\cdot s)=f_{(r)}(m)\cdot s^2$, pour $m\in M$, $r,s\in R_e$.
\end{enumerate}
\end{Definition}
\end{minipage}
\quad 
\begin{minipage}[t]{0.48\textwidth}
 \begin{Definition}\label{df:appCPqu}
 Soit $\RR$ un anneau carr\'e commutatif et soient $(M,A)$ et $(N,B)$ deux
$\RR$-CP-modules, on appelle application $\RR$-CP-\emph{quadratique} de
$(M,A)$ dans $(N,B)$ une application $f:M\to N$ qui satisfait :
\begin{enumerate}
 \item l'image $f(A)$ et celles des d\'efauts : $\im_{\RR} d_f$, $\im_{\RR}
f_{(r)}$ et $\im_{\RR} f_{[x]}$ sont contenues dans $B$, pour $r\in R_e$ et
$x\in R_{ee}$,
 \item le d\'efaut additif $d_f$ et les d\'efauts des commutateurs 
$f_{[x]}$  sont $\RR$-bilin\'eaires,
\item les d\'efauts scalaires $f_{(r)}$ sont homog\`enes de degr\'e 2, i.e.
$f_{(r)}(m\cdot s)=f_{(r)}(m)\cdot s^2$, 
 \item tous les d\'efauts s'annulent sur $A$, c.a.d.
\begin{align*}
d_f(M,A)=d_f(A,M)&=0\\ f_{(r)}(A)&=0\\ f_{[x]}(M,A)=f_{[x]}(A,M)&=0.
\end{align*}
\end{enumerate}
pour $m\in M$, $r,s\in R_e$ et $x\in R_{ee}$.
\end{Definition}

Puisque $f(A)\subset B$, que $d_f$  est biadditif et que
$d_f(M,A)=d_f(A,M)=0$, l'application 
$f\colon (M,A) \to (N,B)$
est une ``quadratic pair map'' au sens de
\cite{QmG}.

\end{minipage}
\medskip


Remarquons qu'une application $\RR$-BHP-quadratique o\`u $\RR$-CP-quadratique annule $0$ car $d_f(0,0)=0$.

\smallskip

\noindent\begin{minipage}[t]{0.48\textwidth}
 \begin{Example}\label{ex:HquadBHP}
 On a vu \eqref{ex:RmodRe} que $R_e$ et $R_{ee}$ sont des $\RR$-BHP-modules
quadratiques. L'appli\-cation $H: R_e\to R_{ee}$ est $\RR$-BHP-quadratique.
On v\'erifie en effet que
 $d_H(r,s)=(s,r)\cdot H(2)$, $H_{[x]}(r,s)=HP((r,s)\cdot x)$, 
$H_{(r)}(s)=(s,s)\cdot H(r)$.
\end{Example}

\end{minipage}
\quad 
\begin{minipage}[t]{0.48\textwidth}
 \begin{Example}\label{ex:HquadCP}
 Soit $A$ un sous-groupe de $R_e$ stable par multiplication \`a droite par
$R_e$ et tel que $P(R_{ee})\subseteq A\subseteq \ZZ_{\RR}(R_e)$, et soit
$B$ un sous-$\ovl{R}$-module de $R_{ee}$. On a vu \eqref{ex:RmodRe} que
$(R_e,A))$ et $(R_{ee},B)$ sont des $\RR$-CP-modules quadratiques. Alors si
$H(A)\subseteq B$ l'application $H: (R_e,A)\to (R_{ee},B)$ est
$\RR$-CP-quadratique.
\end{Example}
\end{minipage}

\begin{Example}\label{ex:muldroite}
Soit $M$ un $\RR$-BHP-module. Pour tout $t\in R_e$ la multiplication \`a
droite $\mu_t: M\to M$, $\mu_t(m)=m\cdot t$ est une application $\RR$-BHP
quadratique. En effet, on a $d_{\mu_t}(m,n)=[m,n]\cdot H(t)$ qui est bien
 $\RR$-bilin\'eaire et d'image $\RR$-centrale. On a aussi
$(\mu_t)_{(r)}(m)= m\cdot(rt-tr)=m\cdot P(y)$ pour un certain $y$. 
L'image est bien $\RR$-centrale et on a $(\mu_t)_{(r)}(m\cdot s)=m\cdot
s\cdot P(y)=   m\cdot P(y)\cdot s^2=(\mu_t)_{(r)}(m)\cdot s^2$. Enfin on a
$(\mu_t)_{[x]}(m,n)= [m,n]\cdot(x\cdot(t-t^2))$ qui  est $\RR$-bilin\'eaire
et d'image $\RR$-centrale.

En particulier, si $M=R_e$, les multiplications \`a droite dans $R_e$
sont des applications $\RR$-BHP-quadratiques.
\end{Example}

La proposition suivante est une cons\'equence imm\'ediate de la d\'efinition \ref{df:appCPqu}.

\begin{Proposition}\label{pr:GrCPquad} Soit $f\colon (M,A) \to (N,B)$ une application $\RR$-CP-quadratique. Elle induit des morphismes de $\ovl{R}$-modules
\[ \ovl{f}={\rm Gr}_1(f) \colon M/A \to N/B\quad\mbox{et}\quad f_2={\rm Gr}_2(f) \colon A \to B\]
tels que $\ovl{f}(\ovl{m}) = \ovl{f(m)}$ et $f_2(a)=f(a)$, pour $m\in M$ et $a\in A$.
\end{Proposition}


\begin{Proposition}\label{pr:factbilin}
 Soient $M$ et $N$ deux $\RR$-BHP-modules, et $\varphi:M\times M\to
\ZZ_{\RR}(N)$ une application $\RR$-bilin\'eaire. Alors $\varphi$ s'annule 
sur $M\times [M]$ et $[M]\times M$. Elle d\'efinit donc une application
$\ovl R$-lin\'eaire $\ovl\varphi: \ovl M\tp_{\ovl R}\ovl M\to
\ZZ_{\RR}(N)$. L'ensemble $\RR\ti Bil(M\times M,\ZZ_{\RR}(N))$ des
applications $\RR$-bilin\'eaires de $M$ dans le $\RR$-centre de $N$ est un
$\overline R$-module.
\end{Proposition}
\begin{preuve}
Soient $m,m',m''\in M$ et $x\in R_{ee}$. Puisque $\varphi$ est
$\RR$-bilin\'eaire, on a  $\varphi([m,m']\cdot x,m'')=
[\varphi(m,m''),\varphi(m',m'')]\cdot x=0$ 
puisque l'image de $\varphi$ est
dans le $\RR$-centre de $N$. Donc l'application injective
$\RR\ti Bil(\ovl{M}\times \ovl{M},\ZZ_{\RR}(N))\to 
\RR\ti Bil(M\times M,\ZZ_{\RR}(N))$ 
induite par la surjection $M\to \ovl{M}$ est bijective. La preuve se
termine gr\^ace \`a la bijection naturelle entre 
$\RR\ti Bil(\ovl{M}\times \ovl{M},\ZZ_{\RR}(N))$ et 
$\ovl R\ti mod(\ovl{M}\otimes_{\ovl{R}} \ovl{M},\ZZ_{\RR}(N))$.
\end{preuve}

\begin{Remark}\label{rm:biaddzen}
  Soit $f:M\to N$ une application $\RR$-BHP-quadratique. Les applications
$d_f$ et $f_{[x]}$ sont $\RR$-bilin\'eaires et
$\RR$-centrales dans l'image. D'apr\`es \ref{pr:factabBHP} elles s'annulent
si l'un de leurs arguments est dans $[M]$, autrement dit
$d_f(M,[M])=d_f([M],M)=0$ et $f_{[x]}(M,[M])=f_{[x]}([M],M)=0$.
\end{Remark}

\begin{Example}\label{ex:commutmul}
 Le commutateur pour la multiplication $R_e\times R_e\to R_e$,
$(r,s)\mapsto rs-sr$ est une application bi-$\RR$-BHP-quadratique. Montrons le pour
la premi\`ere variable: fixons $s$, alors
$c_s(r):=rs-sr=\mu_s(r)-g_s(r)$. Puisque l'anneau $\ovl R$ est
commutatif, l'image de $c_s$ est incluse dans l'image de $P$, donc tous
les  d\'efauts de $c_s$ sont $\RR$-centraux dans l'image de $c_s$. Comme le
groupe additif de $R_e$ est commutatif,
$d_{c_s}=d_{\mu_s}-d_{g_s}=d_{\mu_s}$ puisque $g_s$ est \RR-BP-lin\'eaire.
Donc $d_{c_s}$ est $\RR$-bilin\'eaire. De m\^eme, on a 
\begin{align*}
 (c_s)_{[x]}(r,r')&=c_s([r,r']\cdot x)-[c_s(r),c_s(r')]\cdot x\\
&=(\mu_s)_{[x]}(r,r')-(g_s)_{[x]}(r,r')+[rs,sr']\cdot x +[sr,r's]\cdot x
=(\mu_s)_{[x]}(r,r')+ [r,r']\cdot(x\cdot 2s^2),
\end{align*}
puisque $g_s$ est $\RR$-BHP-lin\'eaire et en utilisant \ref{lm:basrelmod}.
D'apr\`es la proposition pr\'ec\'edente, $c_s$ est $\RR$-bilin\'eaire.
Enfin on a 
\begin{align*}
 (c_s)_{(t)}(r)&=c_s(rt)-c_s(r)t=\mu_s(rt)-g_s(rt)-
(\mu_s(r)t+(-g_s(r))t+ [\mu_s(r),-g_s(r)]\cdot H(t))\\
&=(\mu_s)_{(t)}(r)-(g_s)_{(t)}(r)+[c_s(r),g_s(r)]\cdot
H(t))=(\mu_s)_{(t)}(r),
\end{align*}
puisque $g_s$ est $\RR$-lin\'eaire. Donc $(c_s)_{(t)}$ est quadratique
homog\`ene.
\end{Example}

Explicitons maintenant les diff\'erentes conditions de la d\'efinition
d'une application BHP-quadratique:
\begin{Lemma}\label{lm:defzen}
Soit $f:M\to N$.
\begin{enumerate}
 \item le d\'efaut additif $d_f$ est $\RR$-central dans $\im_{\RR} f$ si et
seulement si
\begin{equation}\label{eq:defaddzen}
 [f(m+m'),f(n)]\cdot x=[f(m),f(n)]\cdot x+[f(m'),f(n)]\cdot x.
\end{equation} 
\item le d\'efaut scalaire $f_{(r)}$ est $\RR$-central dans $\im_{\RR} f$
si et seulement si
\begin{equation}\label{eq:dfscazen}
 [f(m\cdot r),f(n)]\cdot x=([f(m),f(n)]\cdot x)\cdot r.
\end{equation} 
\item le d\'efaut du commutateur $f_{[x]}$ est $\RR$-central dans
$\im_{\RR} f$ si et seulement si
\begin{equation}\label{eq:dfcomzen}
 [f([m,m']\cdot x),f(n)] \cdot y=0.
\end{equation} 
\end{enumerate}
Ces trois relations \'equivalent \`a dire que les applications $m\mapsto
[f(m),n]\cdot x$ sont $\RR$-lin\'eaires pour $n\in \im f$ et $x\in R_{ee}$.
\end{Lemma}
\begin{preuve}
 (1) On doit avoir $[f(m+m')-f(m')-f(m),f(n)]\cdot x=0$. L'additivit\'e du
crochet en la premi\`ere variable (\ref{ax:MC5}) permet de conclure.

(2)  On doit avoir $[f(m\cdot r)-f(m)\cdot r,f(n)]\cdot x=0$. La m\^eme
additivit\'e et (\ref{ax:MC6}) donnent la relation cherch\'ee.

(3) On doit avoir $[f([m,m']\cdot x)-[f(m),f(m')]\cdot x,f(n)]\cdot y=0$.
La m\^eme additivit\'e et (\ref{ax:MC7}) donnent la relation cherch\'ee.
\end{preuve}

\begin{Lemma}\label{lm:dfaddbilin}
 Soit $f:M\to N$ telle que son d\'efaut additif $d_f$ soit $\RR$-central
dans $\im_{\RR} f$.
\begin{enumerate}
 \item L'image du d\'efaut additif est annul\'ee par l'image de $P$.
 \item Le d\'efaut $d_f$ est biadditif si et seulement si
\begin{multline}\label{eq:dfaddbiadd}
 f(m+m'+m'')+f(m)+f(m')+f(m'')=\\
f(m+m')+f(m'+m'')+f(m+m'')-[f(m'),f(m+m'')].
\end{multline}
\item Le d\'efaut $d_f$ est $R_e$-\'equivariant si et seulement si 
\begin{multline}\label{eq:dfaddbieqv}
 f(m\cdot r+n\cdot s)=f(m+n)\cdot rs-f(n)\cdot rs -f(m)\cdot rs \\
+f(m\cdot r) +f(n\cdot s) -[f(m),f(n)]\cdot H(rs).
\end{multline}
\item Le d\'efaut $d_f$ commute aux commutateurs si et seulement si
\begin{equation}\label{eq:dfaddbicomm}
 f([m,m']\cdot x+n)=f([m,m']\cdot x)+f(n)
\end{equation}

\end{enumerate}
\end{Lemma}
\begin{preuve}
(1) En effet, on a d'apr\`es \ref{ax:MC3}
\begin{align*}
 d_f(m,m')\cdot P(x)=[d_f(m,m'),d_f(m,m')]\cdot x=0,
\end{align*}
puisque $d_f(m,m')$ est $\RR$-central.

(2) Pour l'additivit\'e en la premi\`ere variable on doit avoir
\begin{align*}
 d_f(m+m',m'')&=d_f(m,m'')+d_f(m',m'')\\
f(m+m'+m'')&-f(m'')-f(m+m')=d_f(m',m'')+d_f(m,m'')\\
f(m+m'+m'')&=f(m+m')+d_f(m',m'')+f(m'')\\
             &\hskip2cm +f(m+m'')-f(m'')-f(m)\\
f(m+m'+m'')&=f(m+m')+f(m'+m'')-f(m'')\\
            &\hskip2cm -f(m')+f(m'')+f(m+m'')-f(m'')-f(m)\\
f(m+m'+m'')&=f(m+m')+f(m'+m'')-f(m')\\
            &\hskip2cm +[f(m''),f(m')]+f(m+m'')-f(m'')-f(m)\\
f(m+m'+m'')&=f(m+m')+f(m'+m'')-f(m')\\
                 &\hskip2cm +f(m+m'')+[f(m''),f(m')]-f(m'')-f(m)\\
f(m+m'+m'')&=f(m+m')+f(m'+m'')+f(m+m'')+[f(m+m''),f(m')]\\
                &\hskip2cm +[f(m''),f(m')]-f(m')-f(m'')-f(m)\\
f(m+m'+m'')&=f(m+m')+f(m'+m'')+f(m+m'')-\\
             &\hskip2cm [f(m'),f(m+m'')]-f(m'')-f(m')-f(m).
\end{align*}
Un calcul analogue pour la seconde variable aboutit \`a la m\^eme
relation, ce qui montre la deuxi\`eme propri\'et\'e.

(3)  Remarquons tout d'abord que puisque $rs-sr\in \im P$ on a, en
utilisant (1), $d_f(m,m')\cdot rs=d_f(m,m')\cdot sr$. On doit alors avoir 
\begin{align*}
 d_f(m\cdot r,n\cdot s)&=d_f(m,n)\cdot rs
=(f(m+n)-f(n)-f(m))\cdot rs\\
&=f(m+n)\cdot  rs +(-f(n))\cdot rs +(-f(m))\cdot rs +([f(m+n),-f(n)]\cdot
H(rs)\\
        &\hskip3cm+[f(m+n)-f(m)]\cdot H(rs)+[-f(n),-f(m)]\cdot H(rs))\\
&=f(m+n)\cdot  rs -f(n)\cdot rs+[f(n),f(n)]\cdot H(rs) -f(m)\cdot rs 
+[f(m),f(m)]\cdot H(rs) \\
        &\hskip3cm -[f(m),f(n)]\cdot H(rs)-[f(n),f(n)]\cdot
H(rs)-[f(m),f(m)]\cdot H(rs)\\
&=f(m+n)\cdot  rs -f(n)\cdot rs -f(m)\cdot rs
        +f(m\cdot r)+f(n\cdot s)-[f(m),f(n)]\cdot H(rs)
\end{align*}
Ce qui montre la troisi\`eme propri\'et\'e.

(4) On doit avoir 
\begin{equation*}
 d_f([m,m']\cdot x,n)=[d_f(m,n),d_f(m',n)]\cdot x
\end{equation*}
Or le second membre est nul puisque l'image de $d_f$ est $\RR$-centrale. On
d\'eveloppe alors le premier membre:
\begin{equation*}
 d_f([m,m']\cdot x,n)=0=f([m,m']\cdot x+n)-f(n)-(f([m,m']\cdot x)
\end{equation*}
ce qui donne la derni\`ere \'equation.
\end{preuve}

\begin{Lemma}\label{lm:dfComBilin}
 Soit $f:M\to N$ 
\begin{enumerate}
 \item Si le d\'efaut $d_f$ est $\RR$-bilin\'eaire et $\RR$-central dans
$\im_{\RR} f$ alors les d\'efauts des commutateurs sont biadditifs.
\item Si les d\'efauts $f_{(r)}$ sont $\RR$-centraux dans $\im_{\RR} f$
alors les d\'efauts des commutateurs sont $R_e$-bi\'equi\-variants si et
seulement si $f_{(r)}([M])=0$ pour tout $r\in R_e$, c.a.d. si et seulement
si 
\begin{equation}\label{eq:frAnulderiv}
 f([m,n]\cdot x\cdot r)=(f([m,n]\cdot x))\cdot r
\end{equation} 
\item Si les d\'efauts $f_{[x]}$ sont $\RR$-centraux dans $\im_{\RR} f$
alors ils sont $R_{ee}$-\'equivariants.
\end{enumerate}
\end{Lemma}
\begin{preuve} 
(1) En effet, en utilisant successivement \ref{lm:defzen}(1)  et
\ref{lm:dfaddbilin}(4), on a:
\begin{align*}
 f_{[x]}(m,m'+m'')&=f([m,m'+m'']\cdot x)-[f(m),f(m'+m'')]\cdot x\\
&=f([m,m']\cdot x+[m,m'']\cdot x)-[f(m),f(m')+f(m'')]\cdot x\\
&=f([m,m']\cdot x)+f([m,m'']\cdot x)-
([f(m),f(m')]\cdot x+[f(m),f(m'')]\cdot x)\\
&=f([m,m']\cdot x)-[f(m),f(m')]\cdot x+f([m,m'']\cdot x)-[f(m),f(m'')]\cdot
x\\
&=f_{[x]}(m,m')+f_{[x]}(m,m'').
\end{align*}
Ce qui montre la premi\`ere partie du lemme.

(2) On a 
\begin{equation*}
 f_{(r)}([m,m']\cdot x)=f(([m,m']\cdot x)\cdot r)-f([m,m']\cdot x)\cdot r
\end{equation*}
Or puisque $\RR$ est commutatif,
\begin{equation*}
 (([m,m']\cdot x)\cdot r)=[m,m']\cdot (x\cdot r)=[m,m']\cdot ((r,1)\cdot x)
=[m\cdot r,m']\cdot x,
\end{equation*}
donc
\begin{align*}
 f_{(r)}([m,m']\cdot x)&=f([m\cdot r,m']\cdot x)-f([m,m']\cdot x)\cdot r\\
&=f_{[x]}(m\cdot r,m')+[f(m\cdot r),f(m')]\cdot x\\
&\hskip 20mm -(f_{[x]}(m,m')+[f(m),f(m']\cdot x)\cdot r
\end{align*}

Comme $[f(m),f(m']\cdot x$ est central, on a
\begin{equation*}
 -(f_{[x]}(m,m')+[f(m),f(m']\cdot x)\cdot r=
-[f(m),f(m']\cdot (x\cdot r)-f_{[x]}(m,m')\cdot r,
\end{equation*}
et en utilisant la relation \eqref{eq:dfscazen} on obtient:
\begin{equation*}
 f_{(r)}([m,m']\cdot x)=f_{[x]}(m\cdot r,m')-(f_{[x]}(m,m'))\cdot r.
\end{equation*}

(3) Supposons maintenant que les images des $f_{[x]}$ sont centrales dans
l'image de $f$. On a d'une part d'apr\`es \ref{lm:defzen}(3)
\begin{align*}
 f_{[x]}([m,m']\cdot y,n)=f([[m,m']\cdot y,n]\cdot x)-
[f([m,m']\cdot y),f(n)]\cdot x=-[f([m,m']\cdot y),f(n)]\cdot x=0,
\end{align*}
et d'autre part:
\begin{align*}
 [f_{[x]}(m,n),f_{[x]}(m',n)].y&=0.
\qedhere\end{align*}
\end{preuve}

On peut r\'esumer  ces trois lemmes sous la forme:
\begin{Proposition}\label{pr:defquadform}
 Une application $f:M\to N$ est $\RR$-BHP-quadratique si et seulement si
elle v\'erifie les relations:
\begin{subequations}
\begin{align}
&[f(m+m'),f(n)]\cdot x=[f(m),f(n)]\cdot x+[f(m'),f(n)]\cdot x,\\
 &[f(m\cdot r),f(n)]\cdot x=([f(m),f(n)]\cdot x)\cdot r,\\
&[f([m,m']\cdot x),f(n)] \cdot y=0,\\
&f(m+m'+m'')+f(m)+f(m')+f(m'')=\label{eq:dfaddbiadd2}\\
&\hskip1cm f(m+m')+f(m'+m'')+f(m+m'')-[f(m'),f(m+m'')],\notag\\
&f(m\cdot r+n\cdot s)=f(m+n)\cdot rs-f(n)\cdot rs -f(m)\cdot rs \\
&\hskip3cm +f(m\cdot r) +f(n\cdot s) -[f(m),f(n)]\cdot H(rs),\notag\\
&f([m,m']\cdot x+n)=f([m,m']\cdot x)+f(n)\label{eq:dfaddbicomm2}\\
&f(m\cdot s r)-(f(m\cdot s))\cdot r=(f(m\cdot r)-(f(m))\cdot r)\cdot s^2,\\
&f([m,n]\cdot x\cdot r)=(f([m,n]\cdot x))\cdot r,
\end{align}
\end{subequations}
pour tous $m,m',m'',n\in M$, $r,s\in R_e$ et $x,y\in R_{ee}$.
\end{Proposition}
Il r\'esulte de cette proposition et plus particuli\`erement du lemme
\ref{lm:dfComBilin}

\noindent\begin{minipage}[t]{0.48\textwidth}
 \begin{Corollary}\label{cor:BHPquad}
 Une application $f:M \to N$  entre deux $\RR$-BHP-modules est
$\RR$-BHP-quadratique si et seulement si elle satisfait les propri\'et\'es
suivantes:
\begin{enumerate}
 \item les applications $m\mapsto
[f(m),n]\cdot x$ sont $\RR$-lin\'eaires pour $n\in \im f$ et $x\in
R_{ee}$.

\kern 4.5mm
 \item le d\'efaut additif $d_f$ est $\RR$-bilin\'eaire,
\item les d\'efauts scalaires $f_{(r)}$ sont quadratiques homog\`enes, i.e.
$f_{(r)}(m\cdot s)=f_{(r)}(m)\cdot s^2$, pour $m\in M$, $r,s\in R_e$,
\item on a $f_{(r)}([M])=0$.
\end{enumerate}
\end{Corollary}
\end{minipage}
\quad 
\begin{minipage}[t]{0.48\textwidth}
 \begin{Corollary}\label{cor:CPquad}
 Une application $f:(M,A) \to (N,B)$ entre deux $\RR$-CP-modules est
$\RR$-CP-quadratique si et seulement si elle satisfait les propri\'et\'es
suivantes:
\begin{enumerate}
 \item l'image $f(A)$ et les images du  
 d\'efaut additif  et des d\'efauts scalaires,
  $\im_{\RR} d_f$ et
$\im_{\RR} f_{(r)}$,   sont contenues dans $B$, pour
$r\in R_e$.

 \item le d\'efaut additif $d_f$ est $\RR$-bilin\'eaire,
 
\item les d\'efauts scalaires $f_{(r)}$ sont quadratiques homog\`enes, i.e.
$f_{(r)}(m\cdot s)=f_{(r)}(m)\cdot s^2$, pour $m\in M$, $r,s\in R_e$,

\item  le d\'efaut additif et les d\'efauts scalaires  s'annulent sur $A$ :
\[
d_f(M,A)=d_f(A,M)=0=f_{(r)}(A)
\]
pour $m\in M$ et $r,s\in R_e$.
\end{enumerate}
\end{Corollary}
\end{minipage}

\medskip

\noindent On notera que la caract\'erisation d'une application $\RR$-CP-quadratique dans le corollaire \ref{cor:CPquad} n'implique aucune condition portant sur les d\'efauts des commutateurs ; en effet, leur propri\'et\'es souhait\'ees d\'ecoulent des conditions portant sur les autres d\'efauts : $\im_{\RR} f_{[x]}  \subset B$
car $f(A)\subset B \supset [N]$, l'application $f_{[x]}$ annule $A$ en chaque variable car les commutateurs de $M$ resp.\ de $N$ annulent $A$ resp.\ $B$ et donc $f(A)$ en chaque variable, et $f_{[x]}$ est $\RR$-bilin\'eaire d'apr\`es le lemme \ref{lm:dfComBilin}.

\smallskip

\begin{Lemma}\label{lm:3defauts}{\rm [Lemme des 3 d\'efauts]}
 Soit $f:M\to N$ une application v\'erifiant les conditions (1) et (2) de
\ref{df:appBHPqu} , et $r\in R_e$ alors
\begin{equation*}
 d_{f_{(r)}}(m,m')=f_{[H(r)]}(m,m')+d_f(m,m')\cdot(r^2-r).
\end{equation*}
\end{Lemma}

\begin{preuve} 
On a  
\begin{align*}
 d_{f_{(r)}}(m,m')&=f_{(r)}(m+m')-f_{(r)}(m')-f_{(r)}(m)\\
&=f((m+m')\cdot r)-(f(m+m'))\cdot r-f_{(r)}(m')-f_{(r)}(m)\\
&=f(m\cdot r + m'\cdot r +[m,m']\cdot H(r))\\
&\hskip2cm -(f(m)+f(m')
 +d_f(m,m'))\cdot r-f_{(r)}(m')-f_{(r)}(m)\\
&=f(m\cdot r + m'\cdot r)+f([m,m']\cdot H(r))\\
&\hskip2cm 
-(f(m)+f(m'))\cdot r-d_f(m,m')\cdot r-f_{(r)}(m')-f_{(r)}(m)\\
&=d_f(m\cdot r,m'\cdot r)+f(m\cdot r)+f(m'\cdot r)+f([m,m']\cdot H(r))\\
&\hskip2cm -[f(m),f(m')]\cdot H(r)-f(m')\cdot r -f(m)\cdot r\\
&\hskip4cm -d_f(m,m')\cdot r-f_{(r)}(m')-f_{(r)}(m)\\
\end{align*}
puisque $f_{(r)}$ et $f_{[H(r)]}$ sont $\RR$-centraux dans $\im_{\RR} f$:
\begin{align*}
 d_{f_{(r)}}(m,m')&=f_{[H(r)]}(m,m')+d_f(m,m')\cdot r^2+f(m\cdot r)
-f(m)\cdot r\\
&\hskip2cm 
-d_f(m,m')\cdot r-f_{(r)}(m)\\
&=f_{[H(r)]}(m,m')+d_f(m,m')\cdot r^2-d_f(m,m')\cdot r\\
&=f_{[H(r)]}(m,m')+d_f(m,m')\cdot (r^2-r)\qedhere
\end{align*}
\end{preuve}

\noindent\begin{minipage}[t]{0.48\textwidth}
 \begin{Proposition}\label{pr:defscalquadBHP}
 Soit $f:M\to N$ une application $\RR$-BHP-quadratique. Les d\'efauts
scalaires $f_{(r)}:M\to N$ sont $\RR$-BHP-quadratiques.
\end{Proposition}
\begin{preuve} D'apr\`es le lemme pr\'ec\'edent et puisque les d\'efauts
additifs et des commutateurs de $f$ commutent entre eux, les d\'efauts
$d_{(r)}$ sont $R_e$-bilin\'eaires.
\end{preuve}
\end{minipage}
\quad 
\begin{minipage}[t]{0.48\textwidth}
 \begin{Proposition}\label{pr:defscalquadCP}
 Soit $f:(M,A)\to (N,B)$ une application $\RR$-CP-quadratique. Les
d\'efauts scalaires $f_{(r)}:(M,A)\to (N,B)$ sont $\RR$-CP-quadratiques.
\end{Proposition}
\begin{preuve} Puisque $f_{(r)}(A)=0$, les conditions suppl\'ementaires de
la d\'efinition sont imm\'ediatement satisfaites.
\end{preuve}
\end{minipage}

\smallskip

Le corollaire \ref{cor:CPquad} et la proposition \ref{pr:defscalquadCP} permettent de donner une caract\'erisation plus concise d'une application $\RR$-CP-quadratique, comme suit.

 \begin{Proposition}\label{pr:quadCPtensGamma}
 Soient $(M,A)$ et $(N,B)$ deux $\RR$-CP-moduleset $f\colon  \to N$ une application. Alors $f\colon (M,A) \to (N,B)$ est $\RR$-CP-quadratique si et seulement si elle satisfait les propri\'et\'es suivantes.
  \begin{enumerate}
 \item $f(A)\subset B$ ;
 
 \item l'application $d_f\colon M\times M \to N$ se factorise par un morphisme de $\ovl{R}$-modules
 \[ \ovl{d_f} \: \colon \: \ovl{M} \otimes_{\ovl{R}} \ovl{M} \:\longrightarrow\: B\,,\]
  i.e.\ $d_f(m,m') = i_B\ovl{d_f}(\ovl{m} \otimes \ovl{m'})$ o\`u $i_B\colon B \hookrightarrow N$ ;
 
 \item pour $r\in R_e$, l'application $f_{(r)} \colon M \to N$ se factorise par un morphisme de $\ovl{R}$-modules
\[ \ovl{f_{(r)}} \: \colon \: \Gamma^2_{\ovl{R}}(M) \:\longrightarrow\: B\,,\]
i.e.\ $f_{(r)} = i_B  \ovl{f_{(r)}} \gamma_2$ o\`u 
 $\Gamma^2_{\ovl{R}}(M)$ est le terme homog\`ene de degr\'e $2$ de l'alg\`ebre des puissances divis\'ees sur $M$, cf.\ \cite{Rob} ou \cite{GH} o\`u sa structure a \'et\'e d\'etermin\'ee pour tout anneau commmutatif $\ovl{R}$.
 
 \end{enumerate}
 
  \end{Proposition} 

\smallskip

\begin{Proposition}\label{pr:2appRqu}
 Soit $f:(M,A)\to (N,B)$ une application $\RR$-CP-qua\-dra\-tique,
l'application $f:M\to N$ est $\RR$-BHP-quadratique.
\end{Proposition}

\begin{preuve}
 C'est imm\'ediat: puisque $B\subset \ZZ_{\RR}(N)$, il centralise aussi 
$f(M)$.
\end{preuve}

\begin{Proposition}\label{pr:rec2appRqu}Inversement soit $f:M\to N$ une
application $\RR$-BHP-quadratique. L'applica\-tion $f: (M,[M])\to
(\im_{\RR} f,\ZZ_{\RR}(\im_{\RR} f))$ est $\RR$-CP-quadra\-tique.
\end{Proposition}
\begin{preuve}
 Remarquons d'abord que, par \ref{cor:BHPvsCPmod}, $(M,[M])$ et 
$(\im_{\RR} f,\ZZ_{\RR}(\im_{\RR} f))$ sont des $\RR$-CP-modules.
Appliquant les corollaires \ref{cor:BHPquad} et \ref{cor:CPquad} il suffit
de montrer que $d_f(M,[M])=d_f([M],M)=f_{[x]}(M,[M])=f_{[x]}([m],M)=0$, ce
qui d\'ecoule de \ref{rm:biaddzen}, et que
 $f([M])\subseteq \ZZ_{\RR}(\im_{\RR} f)$. Or on a
\begin{equation*}
 [f([m,m']\cdot x),f(m'')]\cdot y=f([[m,m']\cdot x,m'']\cdot y) -
f_{[y]}([m,m']\cdot x,m'').
\end{equation*}
Le premier terme de la diff\'erence est nul par
\ref{ax:MC7}, et le second par la remarque \ref{rm:biaddzen}. 
\end{preuve}

\begin{Remark}\label{rem:defcomm} Soit $f$ une application
$\RR$-BHP-quadratique, on peut v\'erifier que
$d_{f_{(2)}}(m,m')=d_f(m,m')+d_f(m',m)$. On en d\'eduit que
lorsque $r=2$, la formule du lemme \ref{lm:3defauts} se r\'eduit \`a
$d_f(m',m)-d_f(m,m')=f([m',m])-[f(m'),f(m)]$,
et on retrouve la formule (9) de \cite{QmG}.
\end{Remark}
\begin{Example}
 Lorsque la multiplication dans $R_e$ n'est pas commutative, l'application
$R_e\to R_e$, $r\mapsto r^2$ n'est pas toujours quadratique. Prenons par
exemple l'anneau carr\'e $\RR=\bigotimes_R$ (\ref{sec:anncarTenseur}).
On peut v\'erifier que la condition de biadditivit\'e du d\'efaut additif
\eqref{eq:dfaddbiadd2} est satisfaite pour les \'el\'ements $(r,s)$,
$(r',s')$ et $(r'',s'')$ de $R_e$ si et seulement si
$2rr's''+2rs'r''+2sr'r''=0$ c'est \`a dire si et seulement si $R$ est de
caract\'eristique 2. La condition \eqref{eq:dfaddbicomm2} pour les
\'el\'ements $m=m'=(1,0)$, $m''=(r,0)$ et $x=(1,0)$ donne $r^2+r=0$. Par
cons\'equent l'application $R_e\to R_e$, $r\mapsto r^2$ est quadratique si
et seulement si l'anneau $R$ est bool\'een, mais alors $R_e$ est lui
m\^eme au anneau commutatif.
\end{Example}

\subsubsection{Applications quadratiques entre modules sur un anneau
commutatif classique}
\label{sec:quadrclassique}

Soient $R$ un anneau commutatif classique, $\RR=(R\to 0\to R)$ l'anneau
carr\'e associ\'e et $M$ et $N$ deux BHP-modules (cf. \ref{sec:classique})
 Une application $f: M\to N$ est $\RR$-BHP-quadratique si et seulement si
elle est quadratique au sens de \cite{GH}: les conditions portant sur les
d\'efauts des crochets sont triviaux et la $\RR$-bilin\'earit\'e se
r\'eduit \`a la $R$-bilin\'earit\'e.

De m\^eme une application $\RR$-CP-quadratique $f:(M,A)\to (N,B)$ est une
application $R$-quadratique au sens de \cite{GH} v\'erifiant les conditions
(1) et (4) de \ref{df:appCPqu}, c'est \`a dire $R$-quadratique au sens de
\cite{GH}, quadratique dans la cat\'egorie {\bf CP} (\cite{QmG}) et
v\'erifiant la condition (4) de \ref{df:appCPqu}.

\subsubsection{Applications quadratiques entre modules sur l'anneau
$\Znil$}
\label{sec:quadrZnil}
Soient $\Znil=\annc{\Z}$ (\ref{sec:Znil}), $M$ et $N$ deux
$\Znil$-BHP-modules, c'est \`a dire deux groupes nilpotents de classe
$\leq 2$, et $f:M\to N$ une application $\Znil$-BHP-quadratique.  D'apr\`es
\ref{rem:defcomm} on a
\begin{equation*}
 f_{[1]}(m,m')=d_f(m,m')-d_f(m',m),
\end{equation*}
par cons\'equent  les d\'efauts des commutateurs sont biadditifs et
$\Znil$-bilin\'eaires. On v\'erifie m\^eme que 
\begin{equation*}
 f_{[1]}(m,[n,n'])=d_f(m,[n,n'])-d_f([n,n'],m)=0,
\end{equation*}
en effet, puisque $d_f$ est bilin\'eaire et d'image centrale: 
\begin{align*}
 d_f(m,[n,n'])&=d_f(m,n+n'-n-n')\\
&=d_f(m,n)+d_f(m,n')-d_f(m,n)-d_f(m,n')=0.
\end{align*}
On a donc pour tout $r\in \Z$, $f_{[r]}(M,[M])=f_{[r]}([M],M)=0$. Il en
r\'esulte qu'une application $\Znil$-BHP-quadratique $f:M\to N$ est
simplement une application quadratique entre groupes nilpotents de classe $2$.
  En fait, toute application quadratique entre groupes $f\colon G \to H$ se ramm\`ene \`a une application quadratique entre groupes nilpotents de classe $2$, ``\`a un morphisme entre les troisi\`emes termes des suites centrales descendantes pr\`es'', comme suit : notons $\gamma_3(G) = [G,[G,G]]$, $\ovl{G}=G/\gamma_3(G)$, et $\pi_G\colon G \to \ovl{G}$ la projection canonique.

\begin{Proposition}\label{gamma3QG}
Soient $G$ un groupe et $q\colon G \to Q(G)$ son application quadratique universelle \cite{QmG}. Alors la restriction $q_3\colon \gamma_3(G) \to Q(G)$ de $q$ \`a $\gamma_3(G)$ est un isomorphisme sur $\gamma_3(Q(G))$.
\end{Proposition}

\begin{preuve}
L'application $q_3$ est un homomorphisme car $d_q$ annule $[G,G]$ en chaque variable. De plus, pour $a,b,c \in G$, on a $q[a,[b,c]] = [qa,[qb,qc]]$, cf.\  \cite{Niq}, donc $\gamma_3(Q(G)) = q(\gamma_3(G))$. Or $q$ est injective car elle admet  une retraction, \`a savoir, l'application $\ovl{1_G}\colon Q(G) \to G$  obtenue en appliquant la propri\'et\'e universelle de $q$ \`a l'application identique $1_G$ de $G$ qui est bien quadratique.
\end{preuve}

\begin{Corollary}\label{quadNil2gen} Toute application quadratique entre groupes $f\colon G \to H$ induit un diagram commutatif de homomorphismes de groupes dont les lignes sont exactes :
\[\xymatrix{
0 \ar[r] & \gamma_3(G) \ar[d]^-{f_3} \ar[r]^-{q_3} & Q(G) \ar[d]^-{\ovl{f}} \ar[r]^-{Q(\pi_G)} & Q(\ovl{G}) \ar[r] \ar[d]^{\ovl{\ovl{f}}} & 0\\
0\ar[r] & \gamma_3(H) \ar[r] & H \ar[r]^-{\pi_H} & \ovl{H} \ar[r] & 0
}\]
\end{Corollary}

On remarque que $\ovl{\ovl{f}}$ correspond \`a une application quadratique de $\ovl{G}$ dans $\ovl{H}$, donc entre groupes nilpotents de classe $\leq 2$.

\begin{preuve} L'application $\ovl{f}$ est obtenue en appliquant la propri\'et\'e universelle de $q$  \`a $f$. D'apr\`es la proposition \ref{gamma3QG} il ne reste qu'\`a montrer que $\Ker{(Q(\pi_G))}  = \gamma_3(Q(G))$. Rappelons de \cite{QmG} que l'on a une suite exacte naturelle
\[ \xymatrix{
0 \ar[r] & G^{ab} \otimes G^{ab} \ar[r]^-{w_G} & Q(G) \ar[r]^-{\ovl{1_G}} & G \ar[r] & 0 }\,.\]
Elle donne naissance au diagram commutatif suivant dont les lignes sont exactes :
\[\xymatrix{
 & G^{ab} \otimes G^{ab} \hspace{2mm} \ar[r]^-{\pi_{Q(G)}w_G} \ar[d]^{\simeq} &\hspace{2mm} Q(G)/\gamma_3(Q(G)) \hspace{2mm}\ar[r]^-{\ovl{\pi_G 1_G}} \ar[d]^{\ovl{Q(\pi_G)}} & \hspace{2mm}G/\gamma_3(G)\hspace{2mm} \ar[r] \ar@{=}[d]& 0 
\\
0 \ar[r] & \ovl{G}^{ab} \otimes  \ovl{G}^{ab} \ar[r]^{w_{\ovl{G}}} & Q(\ovl{G}) \ar[r]^-{\ovl{1_{\ovl{G}}}} &
\ovl{G}\ar[r] & 0
}\]
Notons que l'application $\ovl{Q(\pi_G)}$ est bien d\'efinie d'apr\`es la proposition \ref{gamma3QG}. On en d\'eduit que $\pi_{Q(G)}w_G$ est injective et que $\ovl{Q(\pi_G)}$ est donc un isomorphisme.
\end{preuve}



\subsection{Module et cat\'egorie des applications CP-quadratiques}
\label{sec:CPmodInt}
\begin{Lemma}\label{lm:CPapp}
 Soit $\RR$ un anneau carr\'e non n\'ecessairement commutatif, soient
$A\subset M$ deux ensembles et $\NB =(N,B)$ un $\RR$-CP-module quadratique.
D\'esignons par $Appl_{A/B}(M,N)$ l'ensemble des applications de $M$ dans
$N$ telles que l'image de $A$ soit contenue dans $B$ et par
$Appl_{A/0}(M,B)$  le sous-ensemble form\'e des applications de $M$ dans
$B$ telles que l'image de $A$ soit 0. D\'efinissons point par point sur ces
ensembles une addition, une action de $R_e$  et une action de $R_{ee}$.
Alors le couple
\begin{equation*}
 Appl(\MA,\NB):=
(Appl_{A/B}(M,N),Appl_{A/0}(M,B))
\end{equation*}
 est un $\RR$-CP-module quadratique.
\end{Lemma}
\begin{preuve}
 Les applications $f+g$, $-f$, $f\cdot r$ et $[f,g]\cdot x$ sont d\'efinies
par 
\begin{align*}
 (f+g)(m)&=f(m)+g(m) & (-f)(m)&=-(f(m)) \\ (f\cdot r)(m)&=(f(m))\cdot r 
& ([f,g]\cdot x)(m)&=[f(m),g(m)]\cdot x\,.
\end{align*}
L'ensemble $Appl_{A/B}(M,N)$ est clairement un sous-groupe du groupe des
applications de $M$ dans $N$, et $Appl_{A/0}(M,B)$ en est un sous-groupe.
Puisque les op\'erations sont d\'efinies point par point elles h\'eritent
de propri\'et\'es du couple $(N,B)$, ce qui montre la proposition.
\end{preuve}

\begin{Proposition}\label{pr:CPmodInt}
 Soient $\MA =(M,A)$ et $\NB =(N,B)$ deux $\RR$-CP-modules
quadratiques, et soit $\ud B:=(B,0)$. Le couple
\begin{equation*}
 \RR\ti\qu(\MA,\NB):=
(\RR\text{\emph{-CP-}quad}(\MA,\NB),\RR\text{\emph{-CP-}quad}(\MA,\ud B))
\end{equation*}
 muni des op\'erations naturelles est un $\RR$-CP-module quadratique.
\end{Proposition}
\begin{preuve}
Montrons que $\RR\text{\emph{-CP-quad}}(\MA,\NB)$ est un sous-groupe  du
groupe des applications de $M$ dans $N$, c.a.d. que $f+g$ et $-f$ sont des
applications $\RR$-CP-quadratiques si $f$ et $g$ le sont. On a 
\begin{align*}
 (f+g)(a)&=f(a)+g(a)\in B & (-f)(a)&=-(f(a))\in B.
\end{align*}
\begin{align*}
 d_{f+g}(m,m')&=(f+g)(m+m')-(f+g)(m')-(f+g)(m)\\
&=f(m+m')+g(m+m')-g(m')-f(m') -g(m)-f(m)\\
&=f(m+m')+d_g(m,m')+[f(m'),g(m)]-f(m')-f(m)\\
&=d_f(m,m')+d_g(m,m')+[f(m'),g(m)]\in B.
\end{align*}
En prenant $g=-f$ on en d\'eduit
\begin{align*}
 d_{-f}(m,m')&=-d_f(m,m')-[f(m),f(m')]\in B.
\end{align*}
\begin{align*}
 (f+g)_{(r)}(m)&=(f+g)(m\cdot r)-((f+g)(m))\cdot r\\
&=f(m\cdot r)+g(m\cdot r)-(f(m)+g(m))\cdot r\\
&=f(m\cdot r)+g(m\cdot r)-g(m)\cdot r -f(m)\cdot r -[f(m),g(m)]\cdot H(r)\\
&=f(m\cdot r)+g_{(r)}(m)-f(m)\cdot r-[f(m),g(m)]\cdot H(r)\\
&=f_{(r)}(m)+g_{(r)}(m)-[f(m),g(m)]\cdot H(r) \in B.
\end{align*}
D'o\`u l'on d\'eduit
\begin{align*}
 (-f)_{(r)}(m)&=-f_{(r)}(m)-[f(m),f(m)]\cdot H(r)\in B.
\end{align*}

Ce qui montre que les applications $f+g$ et $-f$ satisfont les conditions
du (1) du corollaire \ref{cor:CPquad}. On en d\'eduit aussi qu'elles
satisfont les conditions (4) car $f(A) \subset B \subset \mathsf{Z}_{\RR}(N)$.

Les calculs ci-dessus montrent que les d\'efauts additifs  de $f+g$ et $-f$ sont des combinaisons lin\'eaires d'applications
$\RR$-bilin\'eaires de $M$ dans $B$. D'apr\`es la proposition
\ref{pr:factbilin} elles sont aussi $\RR$-bilin\'eaires, ce qui d\'emontre
la condition (2) du corollaire \ref{cor:CPquad}.

Montrons la condition (3). D'apr\`es les calculs ci-dessus, et en utilisant
\ref{pr:defquadform}(b)  on a:
\item \begin{align*}
 (f+g)_{(r)}(m\cdot s)
&=f_{(r)}(m\cdot s)+g_{(r)}(m\cdot s)-[f(m\cdot s),g(m\cdot s)]\cdot H(r)\\
&=f_{(r)}(m)\cdot s^2+g_{(r)}(m)\cdot s^2-([f(m),g(m)]\cdot H(r))\cdot
s^2\\
&=(f_{(r)}(m)+g_{(r)}(m)-[f(m),g(m)]\cdot H(r))\cdot s^2=
((f+g)_{(r)}(m))\cdot s^2\,,
\end{align*}
la factorisation se faisant sans terme correctif puisque tous les termes
sont dans $B$. On a \'egalement:
\begin{align*}
 (-f)_{(r)}(m\cdot s)
&=-f_{(r)}(m\cdot s)-[f(m\cdot s),f(m\cdot s)]\cdot H(r)\\
&=-f_{(r)}(m)\cdot s^2-([f(m),f(m)]\cdot H(r)))\cdot s^2\\
&=((-f)_{(r)}(m))\cdot s^2.
\end{align*}
On a donc montr\'e que
$(\RR\text{\emph{-CP-quad}}(\MA,\NB),\RR\text{\emph{-CP-quad}}(\MA,\ud
B))$ est un
sous-$\RR$-CP-module quadratique de $(Appl_{A/B}(M,N),Appl_{A/0}(M,B))$.
\end{preuve}

\begin{Proposition}\label{pr:composQuad}
 Soient $f:(M,A)\to (N,B)$ et $g:(N,B)\to (L,C)$ deux applications
$\RR$-CP-quadratiques. L'application compos\'ee $g\circ f: (M,A)\to (L,C)$
est $\RR$-CP-quadratique.
\end{Proposition}

\begin{preuve}

D'apr\`es \cite{QmG} le d\'efaut additif de $g\circ f$ est donn\'e par
\begin{equation}
d_{g\circ f}(m,m') = gd_f(m,m') + d_g(f(m),f(m')) \,.\end{equation}
Donc $d_{g\circ f}$ a son image dans $C$, annule $A$ en chaque variable, et est $\RR$-bilin\'eaire d'apr\`es les propositions \ref{pr:GrCPquad} et \ref{pr:quadCPtensGamma}.
 

Consid\'erons les d\'efauts scalaires  
$(g\circ f)_{(r)}$. On a 
\begin{align*}
 (g\circ f)_{(r)}(m)&=(g\circ f)(m\cdot r)-((g\circ f)(m))\cdot r
=g(f(m\cdot r))-(g(f(m))\cdot r\\
&=g(f_{(r)}(m)+f(m)\cdot r)-(g(f(m))\cdot r\\
&=g(f_{(r)}(m))+g(f(m)\cdot r)-(g(f(m))\cdot r+d_g(f_{(r)}(m),f(m)\cdot r)
\end{align*}
Le dernier terme est nul puisque son premier argument est dans $B$, donc
\begin{equation}
 (g\circ f)_{(r)}(m)=g(f_{(r)}(m))+g_{(r)}(f(m)).
\end{equation}
Ce qui montre que $(g\circ f)_{(r)}$ a son image dans $C$ et s'annule si
$m$ est dans $A$. Enfin on a 
\begin{align*}
 (g\circ f)_{(r)}(m\cdot s)&=g(f_{(r)}(m\cdot s))+g_{(r)}(f(m\cdot s))\\
&=g(f_{(r)}(m)\cdot s^2)+g_{(r)}(f_{(s)}(m)+f(m)\cdot s)\\
&=g(f_{(r)}(m))\cdot s^2+g_{(s^2)}(f_{(r)}(m))+g_{(r)}(f_{(s)}(m))\\
&\hskip5cm + g_{(r)}(f(m)\cdot s)+d_{g_{(r)}}(f_{(s)}(m),f(m)\cdot s)
\end{align*}
Or le deuxi\`eme, le troisi\`eme et le cinqui\`eme terme sont nuls puisque
$f_{(r)}(m)$ et $f_{(s)}(m)$ sont dans $B$, donc
\begin{align*}
 (g\circ f)_{(r)}(m\cdot s) =g(f_{(r)}(m))\cdot s^2+g_{(r)}(f(m)\cdot s)\\
&=g(f_{(r)}(m))\cdot s^2+g_{(r)}(f(m)) \cdot s^2=(g\circ f)_{(r)}(m)\cdot
s^2
\end{align*}

Pour des usages futurs, d\'eterminons enfin les d\'efauts $(g\circ f)_{[x]}$. On a
\begin{align*}
 (g\circ f)_{[x]}(m,m')
&=(g\circ f)([m,m']\cdot x)-[(g\circ f)(m),(g\circ f)(m')]\cdot x\\
&=g(f([m,m']\cdot x))-[g(f(m)),g(f(m'))]\cdot x\\
&=g(f_{[x]}(m,m'))+[f(m),f(m')]\cdot x)-[g(f(m)),g(f(m'))]\cdot x\\
&=g(f_{[x]}(m,m'))+g([f(m),f(m')]\cdot x)-[g(f(m)),g(f(m'))]\cdot x\\
&\hskip5cm +d_g(f_{[x]}(m,m'),[f(m),f(m')]\cdot x)\\
\end{align*}
Le dernier terme est nul puisque son second argument est dans $B$, donc
\begin{equation}
 (g\circ f)_{[x]}(m,m') =g(f_{[x]}(m,m'))+g_{[x]}(f(m),f(m'))\qedhere
\end{equation}

\end{preuve}

\begin{Proposition}\label{pr:compCPquad}
 Soient $f:(M,A)\to (N,B)$ et $g:(N,B)\to (L,C)$ deux applications
$\RR$-CP-quadratiques. L'application
\begin{align*}
 f^*:\RR\ti\qu((N,B),(L,C))&\to \RR\ti\qu((M,A),(L,C))& g&\mapsto g\circ f
\end{align*}
est $\RR$-CP-lin\'eaire. Et l'application
\begin{align*}
 g_*:\RR\ti\qu((M,A),(N,B))&\to \RR\ti\qu((M,A),(L,C))& f&\mapsto g\circ f
\end{align*}
est $\RR$-CP-quadratique.
\end{Proposition}
\begin{preuve}
 Montrons que $f^*$ est additive. On a 
\begin{align*}
 (f^*(g_1+g_2))(m)&=((g_1+g_2)\circ f)(m)
=(g_1+g_2)(f(m))=g_1(f(m))+g_2(f(m))\\
&= (g_1\circ f)(m)+(g_2\circ f)(m)=(g_1\circ f+g_2\circ f)(m)
=(f^*(g_1)+f^*(g_2))(m).
\end{align*}
Montrons que $f^*$ est $R_e$-\'equivariante. On a 
\begin{align*}
 (f^*(g\cdot r))(m)&=((g\cdot r)\circ f)(m)=(g\cdot r)(f(m))=g(f(m))\cdot
r 
=(f^*(g))(m)\cdot r=((f^*(g))\cdot r)(m).
\end{align*}
Pour montrer que $f^*$ est $\RR$-lin\'eaire il reste \`a montrer que $f^*$
conserve les crochets. On a 
\begin{align*}
 (f^*([g_1,g_2]\cdot x))(m)&=([g_1,g_2]\cdot x)(f(m))=
[g_1(f(m)),g_2(f(m))]\cdot x\\
&=[(f^*(g_1))(m)),(f^*(g_2))(m)]\cdot x=
([(f^*(g_1)),(f^*(g_2))]\cdot x)(m).
\end{align*}
Montrons maintenant que $g_*$ est $\RR$-CP-quadratique. On a
\begin{align*}
 (d_{g_*}(f_1,
f_2))(m)&=(g_*(f_1+f_2)-g_*(f_2)-g_*(f_1))(m)=(g(f_1+f_2))(m)-
(g(f_2))(m)-(g(f_1))(m)\\
&=g(f_1(m)+f_2(m))-g(f_2(m))-g(f_1(m))\\
&=g(f_1(m))+g(f_2(m))+d_g(f_1(m),f_2(m))-g(f_2(m))-g(f_1(m))\\
=d_g(f_1(m),f_2(m))
\end{align*}
Il en r\'esulte imm\'ediatement que $d_{g_*}$ est $\RR$-bilin\'eaire en les
variables $f_1$ et $f_2$, qu'elle s'annule si l'image de $f_1$ ou celle de
$f_2$ est dans $B$, que son image est constitu\'ee d'applications de $M$ dans $C$. Il reste \`a
\'etudier les d\'efauts scalaires.
\begin{align*}
 (g_*)_{(r)}(f)&=g_*(f\cdot r)-(g_*(f))\cdot r=g\circ(f\cdot r)-(g\circ
f)\cdot r.
\end{align*}
\begin{align*}
 (g_*)_{(r)}(f)(m)&=g((f\cdot r)(m))-((g\circ f)\cdot r)(m)
=g(f(m)\cdot r)-(g(f(m))\cdot r=g_{(r)}(f(m))
\end{align*}
Ce qui assure que $(g_*)_{(r)}$ s'annule si l'image de $f$ est dans $B$ et
que l'image de $(g_*)_{(r)}(f)$ est dans $C$. Enfin on a 
\begin{align*}
 (g_*)_{(r)}(f)(m\cdot s)&=g_{(r)}(f(m\cdot s))=g_{(r)}(f(m)\cdot
s+f_{(s)}(m))\\
&=g_{(r)}(f(m)\cdot s)+g_{(r)}(f_{(s)}(m))+
d_{g_{(r)}}(f(m)\cdot s,f_{(s)}(m)).
\end{align*}
Dans cette somme on remarque que $f_{(s)}(m)\in B$ d'apr\`es
\ref{df:appCPqu} (1). Il en r\'esulte d'une part que le deuxi\`eme terme
est nul d'apr\`es \ref{df:appCPqu} (4), et que le dernier terme est aussi
nul car $g_{(r)}$ est quadratique d'apr\`es \ref{pr:defscalquadCP}. On a
donc
\begin{equation}
(g_*)_{(r)}(f)(m\cdot s) =g_{(r)}(f(m)\cdot s)
\end{equation}
En prenant $s=1$ on obtient pour $t\in R_e$
\begin{align*}
(g_*)_{(r)}(f\cdot t)(m) &= g_{(r)}((f\cdot t)(m)) = g_{(r)}(f(m)\cdot t) = 
g_{(r)}(f(m))\cdot t^2 = (g_*)_{(r)}(f)(m)\cdot t^2\\
&=(((g_*)_{(r)}(f))\cdot t^2 )(m)
\end{align*}
On a donc bien $(g_*)_{(r)}(f\cdot t)=((g_*)_{(r)}(f))\cdot t^2$.
\end{preuve}

Les r\'esultats de ce paragraphe se r\'esument comme suit :

\begin{Theorem}\label{cor:CPR}
Il existe une cat\'egorie
$\cpr$  dont les objets sont les $\RR$-CP-modules
quadratiques et dont les morphismes sont les applications
$\RR$-CP-quadratiques. Elle est pr\'equadratique \`a droite au sens de \cite{BHP}, et admet un
foncteur Hom interne.
\end{Theorem}


\begin{Remark}
Soit $R$ un anneau commutatif. En prenant $\RR=(R\to 0\to R)$, on obtient une cat\'egorie $\mathsf{CP}_R = \cpr$ qui contient toute application $R$-quadratique $f\colon M \to N$ entre $R$-modules $M$, $N$ (au sens de \cite{GH}) comme morphisme, d'au moins deux fa\c{c}ons canoniques :  $f\colon (M,0)\to (N,N)$ et $f\colon (M,rad_R(f)) \to (N,f(rad(f)) + D_f)$ sont des morphismes dans $\mathsf{CP}_R$ o\`u $rad_R(f)$ d\'esigne l'ensemble des \'el\'ements $m$ de $M$ tels que $d_f(m,m') =0$ pour tout $m'\in M$ et $f(mr) =f(m)r$ pour tout $r\in R$, 
et o\`u $D_f$ d\'esigne le sous-$R$-module de $N$ engendr\'e par les images de $d_f$ et des $f_{(r)}$, $r\in R$. Nous nous attendons \`a ce que l'\'etude des cat\'egories $\cpr$ fournisse  alors  une approche nouvelle des applications quadratiques entre modules (carr\'es et classiques) qui sera d\'evelopp\'ee dans \cite{B2}. 
\end{Remark}


\end{document}

\noindent\begin{minipage}[t]{0.48\textwidth}
 
\end{minipage}
\quad 
\begin{minipage}[t]{0.48\textwidth}
 
\end{minipage}
\smallskip